\documentclass[12pt, a4paper]{article}

\usepackage[english]{babel} % English language hyphenation

\usepackage{microtype} % Better typography

\usepackage{amsmath,amsfonts,amsthm} % Math packages for equations

\usepackage[svgnames]{xcolor} % Enabling colors by their 'svgnames'

\usepackage[hang, small, labelfont=bf, up, textfont=it]{caption} % Custom captions under/above tables and figures

\usepackage{booktabs} % Horizontal rules in tables

\usepackage{lastpage} % Used to determine the number of pages in the document (for "Page X of Total")

\usepackage{graphicx} % Required for adding images

\usepackage{enumitem} % Required for customising lists
\setlist{noitemsep} % Remove spacing between bullet/numbered list elements

\usepackage{sectsty} % Enables custom section titles
\allsectionsfont{\usefont{OT1}{phv}{b}{n}} % Change the font of all section commands (Helvetica)

\usepackage{geometry} % Required for adjusting page dimensions

\usepackage[T1]{fontenc} % Output font encoding for international characters
\usepackage[utf8]{inputenc} % Required for inputting international characters

\usepackage{XCharter} % Use the XCharter font

\usepackage{fancyhdr} % Needed to define custom headers/footers
\fancypagestyle{firstpage}{ % Page style for the first page with the title
	\fancyhf{}
}

\usepackage{titling} % Allows custom title configuration

\newcommand{\HorRule}{\color{DarkGoldenrod}\rule{\linewidth}{1pt}} % Defines the gold horizontal rule around the title

\pretitle{
	\vspace{-30pt} % Move the entire title section up
	\HorRule\vspace{10pt} % Horizontal rule before the title
	\fontsize{32}{36}\usefont{OT1}{phv}{b}{n}\selectfont % Helvetica
	\color{DarkRed} % Text colour for the title and author(s)
}

\posttitle{\par\vskip 15pt} % Whitespace under the title

\preauthor{} % Anything that will appear before \author is printed

\postauthor{ % Anything that will appear after \author is printed
	\vspace{10pt} % Space before the rule
	\par\HorRule % Horizontal rule after the title
	\vspace{20pt} % Space after the title section
}

\usepackage{lettrine} % Package to accentuate the first letter of the text (lettrine)
\usepackage{fix-cm}	% Fixes the height of the lettrine

\newcommand{\initial}[1]{ % Defines the command and style for the lettrine
	\lettrine[lines=3,findent=4pt,nindent=0pt]{% Lettrine takes up 3 lines, the text to the right of it is indented 4pt and further indenting of lines 2+ is stopped
		\color{DarkGoldenrod}% Lettrine colour
		{#1}% The letter
	}{}%
}

\usepackage{xstring} % Required for string manipulation

\newcommand{\lettrineabstract}[1]{
	\StrLeft{#1}{1}[\firstletter] % Capture the first letter of the abstract for the lettrine
	\initial{\firstletter}\textbf{\StrGobbleLeft{#1}{1}} % Print the abstract with the first letter as a lettrine and the rest in bold
}

\usepackage{amssymb}

\newtheorem{theorem}{Theorem}  
\newtheorem{proposition}{Proposition}
\newtheorem{corollary}{Corollary}
\newtheorem{lemma}{Lemma}

\usepackage{amsbsy}

\usepackage{tikz-cd}

\title{On the geometry and topology of Cayley condition in Poncelet porism for triangles}

\author{{\bf Yirmeyahu Kaminski}\\
		School of Mathematical Sciences,\\
		Holon Institute of Technology, Holon, Israel. \\
		email: kaminsj@hit.ac.il}

\date{}

\begin{document}

%----------------------------------------------------------------------------------------

\maketitle % Print the title

\thispagestyle{firstpage} % Apply the page style for the first page (no headers and footers)

%----------------------------------------------------------------------------------------
%	ABSTRACT
%----------------------------------------------------------------------------------------

\lettrineabstract{}

\begin{abstract}
	This paper investigates the differential-geometric and topological
	properties of the Cayley condition in Poncelet porism for triangles, defined
	as the locus of pairs of non-degenerate conics that admit a Poncelet triangle.
	While the algebraic condition for this porism, established by Cayley, is
	classical, the geometric nature of the set of solutions has remained largely
	unexplored. We demonstrate that this Cayley set is a smooth, connected,
	9-dimensional complex manifold. This is proven by showing it is an open
	subset of a smooth algebraic variety endowed with a trivial fiber bundle
	structure over the space of non-degenerate conics.
	
	To further analyze its structure, we construct the moduli space of transversely
	intersecting conic pairs under the action of
	$\mathbb{P}GL_{3}(\mathbb{C})$ and identify it with an open subset of
	$\mathbb{CP}^{2}/S_{3}$. We compute the fundamental group of a generic
	orbit.
	
	The elliptic j-invariant is then introduced as a holomorphic map on this
	space of conics, which factors through this moduli space. We analyze the
	subset of the Cayley set where this map is a submersion---its regular part---which
	corresponds to excluding points whose j-invariant is one of the critical
	values $0$ or $1728$. We prove that this regular part is the total space
	of a fiber bundle over $\mathbb{C}\setminus \{0,1728\}$. This structure
	allows for the computation of the fundamental group via the long exact
	sequence of homotopy.
	
	Finally, we provide a principal bundle formulation for the Poncelet correspondence
	itself over orbits of conic pairs.
\end{abstract}

\bigskip
\textbf{MSC:} 51N35, 14H52, 32Q55, 55R10, 57M05

%%*************** Text entry area ******************%%
\section{Introduction}
\label{sec1}

Poncelet's porism, dating back to the early 19th century, represents one
of the most elegant results in projective geometry. The theorem states
that given two non-degenerate conics, if there exists one n-sided polygon
inscribed in one conic and circumscribed about the other, then there exist
infinitely many such polygons. This remarkable closure property has attracted
significant attention from mathematicians for over two centuries and continues
to reveal deep connections between algebraic geometry, differential geometry,
and complex analysis. A thorough treatment with many extensions appears
in~\cite{Dragovic-Radnovic-2011}.

Cayley's contribution, formulated in the mid-19th century, provided an
algebraic characterization of when such polygons exist. Through the polynomial
expansion of $\sqrt{\det (tC + D)}$, where $C$ and $D$ are the matrices
of the conics, Cayley established precise conditions on the coefficients
that determine the existence of n-sided Poncelet polygons. For triangles,
this condition is particularly elegant: the coefficient of $t^{2}$ in the
Taylor expansion must vanish.

Despite the extensive literature on Poncelet's porism and Cayley's conditions,
several fundamental geometric questions remain unexplored. In particular,
the differential-geometric and topological properties of the set of conic
pairs satisfying Cayley's conditions have not been thoroughly investigated.
This paper addresses this gap by examining the smoothness and topological
properties of the Cayley set---the locus of conic pairs that admit Poncelet
polygons. In this paper, we focus on the case of triangles.
There is however, one major exception: in~\cite{Jakob-93},
a fine moduli space is constructed for such pairs of conics with a fixed
bitangent. Given the sparsity of the literature on the subject, we will
devote section~\ref{sec::previous_work} to a more detailed review of~\cite{Jakob-93}
and a comparison with our work.

Our main contributions are as follows. First in section~\ref{sec::manifold},
we provide an explicit formulation of the Cayley equation for triangles
in terms of the coefficients of the conics. Second, we prove that the Cayley
set is a complex manifold by showing it is an open subset of a smooth complex
algebraic variety with a natural fiber bundle structure. We also show that
it is connected.

Then in section~\ref{sec::pair_of_conics}, we turn our attention to the
space of pairs of non-degenerate conics that intersect transversely. We
show that the group $\mathbb PGL_{3}(\mathbb C)$ naturally acts on it.
We compute the isotropy groups associated to this action. We compute the
fundamental group of a generic orbit. We show how the space of orbits can
be identified to an open set of the quotient of the projective space by
the symmetry group $S_{3}$, $\mathbb C\mathbb P^{2}/S_{3}$. We show that
the projection is a holomorphic submersion.

In section~\ref{sec::elliptic_curves} we then exhibit the relation between
this construction and elliptic curves and show that the j-invariant defines
a submersion from the above space of pairs of conics, outside a critical
set.

Section~\ref{sec::cayley_set_j_inv} is devoted to show the action of
$\mathbb PGL_{3}(\mathbb C)$ is well defined on the Cayley set and that
the j-invariant yields a fiber bundle structure on the Cayley set outside
some critical set defined as the inverse image of the critical values
$0$ and $1728$ of the j-invariant. This allows to extract some information
about the fundamental group of the set of regular points.

Finally in section~\ref{sec::poncelet_principal_bundle}, we introduce a
principal bundle formulation of the Poncelet correspondence.
	
\subsection{Previous work on moduli space of Poncelet pairs}
\label{sec::previous_work}

In this section, we shall give more detail about~\cite{Jakob-93}. The purpose
of this paper was to construct a fine moduli space of Poncelet pairs of
level $n$ with fixed bitangent for $n \geq 4$. Here a pair of conics
$Q$ and $Q'$ in the complex projective space is called a
\textit{Poncelet pair of level $n$}, if the conics intersect transversely
and if there is a Poncelet polygon, defined by a sequence of lines and
points $(l_{1},p_{1,2}), \ldots , (l_{n},p_{n,1})$ such that the lines
$l_{i}$ are tangent to $Q$, the points $p_{i,j}$ lie in $Q'$ and
$p_{i,j} = l_{i} \cap l_{j}$. It is proven that the quotient of the upper
half plane $\mathcal{H} \subset \mathbb C$ by the group
\begin{equation*}
	\Gamma _{1} = \left \lbrace \left (
	\begin{matrix}
		a & b
		\\
		c & d
	\end{matrix}
	\right ) \in SL(2,\mathbb Z) \, | \,c \equiv 0 [n], a \equiv d
	\equiv 1 [n] \right \rbrace ,
\end{equation*}
namely $Y_{1}(n) = \mathcal{H}/\Gamma _{1}(n)$ is a fine moduli space for
Poncelet pairs with fixed bitangent, for $n \geq 4$.

The starting point of the proof is the existence of a universal family
$\tilde{Y}_{1}(n) \rightarrow Y_{1}(n)$ of elliptic curves with a point
of order $n$, such that every flat family (see~\cite{Vakil-25}) of elliptic
curves with a point of order $n$ over a variety $S$ is a pullback of this
universal family under a uniquely determined morphism
$S \rightarrow Y_{1}(n)$. Then $\tilde{Y}_{1}(n)$ is proven to be isomorphic
to a universal family of Poncelet pairs of level $n$. The whole construction
is algebraic and does not provide any analysis of the topological and differential
properties of the moduli space. Moreover this approach is valid for
$n \geq 4$ and requires the choice of a bitangent.

By contrast, the present work does not introduce auxiliary marked data
(a bitangent, a base point, a non-halfperiod point) and focuses directly
on the total space of pairs of conics $(C, D)$ that intersect transversely
and on the Cayley locus for $n=3$ (which is excluded from~\cite{Jakob-93}),
viewed as a complex manifold in the classical topology, and we describe
its differential-topological structure: smoothness, connectedness, the
construction of fiber bundles and a holomorphic submersion and a covering
space involving the Cayley set, as detailed below. We also compute the
fundamental group of some sets related to it. None of this appears in~\cite{Jakob-93}.
However the use of elliptic curves is of course omnipresent in both approaches.

It is worth mentioning~\cite{BarthMichel-93}, where the problem is approached
from the perspective of enumerative geometry, utilizing rational elliptic
surfaces and modular curves to count the number of conics in a general
pencil that are $n$-inscribed or $m$-circumscribed. This of course is far
from the scope and methods of our work.
	
\section{Cayley's theorem}
\label{sec::cayley_theorem}
	
First let us recall what is meant by Poncelet's porism. As mentioned in~\cite{griffiths_harris_78},
we consider two non-degenerate conics $C$ and $D$ in the complex projective
plane. These two conics are assumed to intersect transversely (as in~\cite[p.38]{griffiths_harris_78}),
so that they meet at four different points. Poncelet's porism states that
either there are no $n$-sided polygons simultaneously inscribed in
$C$ and circumscribed about $D$, or there are infinitely many of them.

Cayley's theorem gives a necessary and sufficient algebraic condition for
the existence of at least one such $n$-sided polygon. See~\cite{griffiths_harris_78}
for more details and a proof of the theorem. A more recent presentation
also appears in~\cite[Section 2.2.1]{Dolgachev-12}. Here and
throughout the paper, we shall identify a conic to its defining matrix
in the standard projective frame of $\mathbb C\mathbb P^{2}$.

\begin{theorem}
\label{thm1}
Let $C$ and $D$ be two non-degenerate conics of
$\mathbb C\mathbb P^{2}$ meeting at four different points. Consider the
following Taylor expansion at $t=0$:

\begin{equation*}
	\sqrt{\det (tC+D)} = A_{0} + A_{1} t + A_{2} t^{2} + \cdots = \sum _{k=0}^{
		\infty }A_{k} t^{k}.
\end{equation*}

There exists a $n$-sided polygon inscribed in $C$ and circumscribed about
$D$ if, and only if,

\begin{equation*}
	\left |
	\begin{array}{c@{\quad}c@{\quad}c}
		A_{2} & \cdots & A_{m+1}
		\\
		\vdots & & \vdots
		\\
		A_{m+1} & \cdots & A_{2m}
	\end{array}
	\right | = 0, \text{ when } n \text{ is odd and } n=2m+1,
	\text{ for some } m \geq 1,
\end{equation*}
\begin{equation*}
	\left |
	\begin{array}{c@{\quad}c@{\quad}c}
		A_{3} & \cdots & A_{m+1}
		\\
		\vdots & & \vdots
		\\
		A_{m+1} & \cdots & A_{2m-1}
	\end{array}
	\right | = 0, \text{ when } n \text{ is even and } n=2m,
	\text{ for some } m \geq 2.
\end{equation*}
\end{theorem}

Following this theorem, a pair of non-degenerate conics in the complex
plane, in a generic configuration, will have a triangle inscribed in
$C$ and circumscribed about $D$ if, and only if, $A_{2} = 0$. In the sequel,
the set defined by this equation is called \textbf{the Cayley set}. Notice
that in the case of the general Cayley set for $n$-sided polygons, Gerbaldi
in~\cite{Gerbaldi-19} could count the number of conics in a given pencil,
that have an inscribed $n$-sided polygon which is circumscribed about a
given conic of the pencil. This result is reviewed in~\cite{Dolgachev-12,Barth-Bauer-96}.
Here we shed some light on some geometric and topological properties of
the Cayley set for triangles.
	
\section{The Cayley set is a manifold}
\label{sec::manifold}
	
The set of conics in $\mathbb C\mathbb P^{2}$ is naturally seen as
$\mathbb C\mathbb P^{5}$. Indeed $3\times 3$ symmetric matrices defined
modulo multiplication by a non-zero scalar are in one-to-one correspondence
with points in $\mathbb C\mathbb P^{5}$. Let $\mathfrak{U}$ be the dense
open set of $\mathbb C\mathbb P^{5}$ of non-degenerate conics. More precisely
if $[x_{0}:x_{1}:x_{2}:x_{3}:x_{4}:x_{5}]$ are the homogeneous coordinates
in $\mathbb C\mathbb P^{5}$, this defines a conic the matrix of which is:
\begin{equation*}
	C = \left (
	\begin{array}{c@{\quad}c@{\quad}c}
		x_{0} & x_{1} & x_{2}
		\\
		x_{1} & x_{3} & x_{4}
		\\
		x_{2} & x_{4} & x_{5}
	\end{array}
	\right ).
\end{equation*}
By definition of non-degeneracy, the conic lies in $\mathfrak{U}$ if
$\det (C) \neq 0$, which shows that $\mathfrak{U}$ is a Zariski open set
and then a dense open set in the classical topology.

In $\mathfrak{U} \times \mathfrak{U}$, let us consider the set of pairs
$(C,D)$ such that $C$ and $D$ intersect in four different points.

For that purpose let $P(t)$ be the polynomial $P(t) = \det (tC+D)$, where
we have, as above, identified a conic with its matrix in the standard
projective frame of
$\mathbb C\mathbb P^{2}$. Then
$P(t) = \sigma _{3,0} t^{3} + \sigma _{2,1} t^{2} + \sigma _{1,2} t +
\sigma _{0,3}$, where

\begin{equation}
	\label{eq::sigma}
	\begin{aligned}
		\sigma _{3,0} & =  \det (C)
		\\
		\sigma _{2,1} & =  \det (c1,c2,d3) + \det (c1,d2,c3) + \det (d1,c2,c3)
		\\
		\sigma _{1,2} & =  \det (c1,d2,d3) + \det (d1,c2,d3) + \det (d1,d2,c3)
		\\
		\sigma _{0,3} & =  \det (D)
	\end{aligned}
	,
\end{equation}
with $c_{i},d_{j}$ being the columns of respectively $C$ and $D$. By definition
the discriminant~\cite{Lang-02,Hassett-07} of $P(t)$, that we shall denote
$\Delta (P)$, is given by:
$\Delta (P) = (r_{1} - r_{2})^{2}(r_{1} - r_{3})^{2}(r_{2} - r_{3})^{2}$,
where $r_{1},r_{2},r_{3}$ are the roots of $P(t)$. From this expression,
it is clear that $\Delta (P) \neq 0$ if, and only if $P$ has distinct roots.
Moreover, relying on the relations between the coefficients and the roots
of a polynomial, we get the following expression:

\begin{equation}
	\label{eq::discriminant}
	\sigma_{3,0}^4 \Delta (P) = \sigma _{2,1}^{2} \sigma _{1,2}^{2} - 4 \sigma _{3,0}
	\sigma _{1,2}^{3} - 4\sigma _{2,1}^{3} \sigma _{0,3} - 27 \sigma _{3,0}^{2}
	\sigma _{0,3}^{2} + 18 \sigma _{3,0}\sigma _{2,1}\sigma _{1,2}\sigma _{0,3}.
\end{equation}
	
\begin{lemma}
	\label{lem::discriminant}
	The set of pairs $(C,D) \in \mathfrak{U} \times \mathfrak{U}$ that intersect
	transversely or equivalently in four different points is an open set
	$\mathfrak{O}$ of $\mathfrak{U} \times \mathfrak{U}$. More precisely it
	is the complement in $\mathfrak{U} \times \mathfrak{U}$ of the hypersurface
	defined by the discriminant of $P(t)$:
	$\mathfrak{O} = \mathfrak{U} \times \mathfrak{U} \setminus Z(\Delta (P))$.
\end{lemma}
\begin{proof}
	By Bezout's theorem, two conics in the complex projective plane always
	meet in four points counted with multiplicities. A point of intersection
	has multiplicity greater than $1$ if the two conics are tangent at this
	point, which yields an algebraic condition of the conics. Therefore the
	pairs of conics that intersect transversely are the points of an open Zariski
	set in $\mathfrak{U} \times \mathfrak{U}$.
	
	More explicitly if we consider the polynomial $P(t)$, the two conics intersect
	transversely if, and only if their intersection is made of four distinct
	points, which is equivalent to say that there are exactly three distinct
	degenerate conics in the pencil defined by $C$ and $D$, which are the three
	pairs of lines passing through the four intersection points. This is equivalent
	to say that the three roots of $P(t)$ are all distinct, or in other words
	that the discriminant of $P(t)$ does not vanish.
	
	Furthermore, the right-hand side of expression~\eqref{eq::discriminant} 
	is a bi-homogeneous polynomial of degree $(6,6)$ in the entries of
	$C$ and $D$. It therefore defines an algebraic variety of
	$\mathfrak{U} \times \mathfrak{U}$, the complement of which is precisely
	the set of pairs of conics that intersect transversely.
\end{proof}
	
This shows that $\mathfrak{O}$ is a Zariski open set of
$\mathfrak{U} \times \mathfrak{U} \subset \mathbb C\mathbb P^{5}
\times \mathbb C\mathbb P^{5}$, and since the product of irreducible varieties
is irreducible, $\mathfrak{O}$ is a dense proper open subset of it, which
makes it a complex manifold in the classical topology, of dimension
$10$.

Hence the Cayley set is defined by
$\mathfrak{C} = \{(C,D) \in \mathfrak{O} | A_{2} = 0\}$.

Let us compute the explicit expression of $A_{2}$.

If $\phi (t) = \sqrt{P(t)}$, we have
$A_{2} = \frac{1}{2} \phi ''(0) = -\frac{1}{8 \sigma _{0,3}^{3/2}}
\sigma _{1,2}^{2} + \frac{1}{2 \sqrt{\sigma _{0,3}}} \sigma _{2,1}$.

Recall that a constructible set is a finite union of locally closed subsets
(see~\cite[p.39]{harris-92}). Here we consider constructible sets in the
Zariski topology.

\begin{proposition}
\label{prop1}
The Cayley set $\mathfrak{C}$ is a constructible set, and more precisely
an algebraic variety in $\mathfrak{O}$ defined by

\begin{equation}
	\label{eq::cayley_equation}
	-\sigma _{1,2}^{2} + 4 \sigma _{0,3} \sigma _{2,1} = 0.
\end{equation}

Therefore the dimension of $\mathfrak{C}$ is 9.
\end{proposition}

\begin{proof}
	For $(C,D) \in \mathfrak{O} $, the coefficient $A_{2}$ vanishes if, and only
	if, we have:
	$8 \sigma _{0,3}^{3/2} A_{2} = -\sigma _{1,2}^{2} + 4 \sigma _{0,3}
	\sigma _{2,1} = 0$, since $\sigma _{0,3} = \det (D) \neq 0$. This a bi-homogeneous
	equation in $C$ and $D$ of degrees $(2,4)$. Therefore $\mathfrak{C}$ is
	indeed an algebraic subvariety of $\mathfrak{O}$. Since
	$\mathfrak{C}$ is defined by a single equation in
	$\mathfrak{O} \subset \mathbb C\mathbb P^{5} \times \mathbb C
	\mathbb P^{5}$, it has dimension 9.
\end{proof}

In the sequel, equation~\eqref{eq::cayley_equation} will be called
\textit{the Cayley equation} and we shall write it $\gamma (C,D) = 0$ for
short.
	
Now we want to prove the smoothness of the Cayley set $\mathfrak{C}$. For
this purpose, we shall consider the superset
$\widetilde{\mathfrak{C}} = \{(C,D) \in \mathbb C\mathbb P^{5}
\times \mathfrak{U} | \gamma (C,D) = 0 \}$.

%t2 #&#
\begin{theorem}
	\label{thm2}
	The set $\widetilde{\mathfrak{C}}$ is the total space of a smooth trivial
	fiber bundle, with a connected general fiber.
\end{theorem}
\begin{proof}
	Let consider the projection
	$\pi : \widetilde{\mathfrak{C}} \rightarrow \mathfrak{U}, (C,D)
	\mapsto D$. The above analysis simply means that the fiber over
	$D \in \text{im}(\pi )$ is a quadric hypersurface of
	$\mathbb C\mathbb P^{5}$. Since in our setting, the base field is
	$\mathbb C$, the map $\pi $ is obviously surjective.
	
	Now, for a given $D \in \mathfrak{U}$ consider the projectivity of
	$\mathbb C\mathbb P^{5}$ defined in such a way that a point in
	$\mathbb C\mathbb P^{5}$, say $C$, viewed as a conic in
	$\mathbb C\mathbb P^{2}$, is mapped to the matrix which columns are
	$D^{-1}c1, D^{-1}c2, D^{-1}c3$. By this isomorphism, that we denote
	$\psi _{D}$, the fiber of $\pi $ over $D$ is mapped to the fiber over
	$I_{3}$, which has equation
	$- c_{11}^{2} - c_{22}^{2} - c_{33}^{2} + 2 (c_{11} c_{22} + c_{11} c_{33}
	+ c_{22} c_{33}) - 4 (c_{12}^{2} + c_{13}^{2} + c_{23}^{2}) = 0$, and which
	is easily seen to be non-singular, and irreducible in the Zariski topology,
	since a smooth quadric is always irreducible, provided the ambient space
	has dimension at least $2$.
	
	Therefore $\pi $ is a smooth surjective map, which fibers are all isomorphic
	to a smooth quadric hypersurface of $\mathbb C\mathbb P^{5}$. Let
	$Q$ be the fiber over $I_{3}$. Then the following isomorphism holds
	$\widetilde{\mathfrak{C}} \simeq \mathfrak{U} \times Q, (C,D)
	\mapsto (D,\psi _{D}(C))$. This map is actually algebraic and invertible,
	the inverse of which is $(D,E) \mapsto (F,D)$ where $F$ is the matrix the
	columns of which are $D e_{1}, D e_{2}, D e_{3}$, with
	$e_{1},e_{2},e_{3}$ being the columns of $E$, viewed as the matrix of a
	conic in $\mathbb C\mathbb P^{2}$. The inverse map is also a morphism of
	algebraic varieties, so that $\widetilde{\mathfrak{C}}$ and
	$\mathfrak{U} \times Q$ are indeed isomorphic varieties.
	
	Since $\mathfrak{U} \times Q$ is smooth, so is
	$\widetilde{\mathfrak{C}}$. Moreover, as claimed,
	$\widetilde{\mathfrak{C}}$ is the total space of a trivial fiber bundle
	over $\mathfrak{U}$, with $Q$ as a general fiber, which is irreducible
	in the Zariski topology and therefore connected in the classical topology.
\end{proof}
	
Then we can deduce the following conclusion on $\mathfrak{C}$.
	
\begin{corollary}
	%%LEAP%%%\label{cor1}
	\label{cor::Caley_set_smoothness}
	The set $\mathfrak{C}$ is a connected complex manifold of dimension 9.
\end{corollary}
\begin{proof}
	Since $\mathfrak{C}$ is an open set of $\widetilde{\mathfrak{C}}$ in the
	classical topology, it is indeed a complex manifold of dimension $9$.
	
	It is connected since it is the complement of an algebraic hypersurface
	in a complex manifold, which means that when considered as a real manifold,
	it is the complement of a codimension $2$ submanifold. One has just to
	apply the lemma below with $M = \mathfrak{C}$ and $N = Z(\det (C))$.
\end{proof}

%l2 #&#
\begin{lemma}
	%%LEAP%%%\label{lem2}
	\label{lem::connectivity_codim2}
	Consider a real smooth manifold $M$ and a submanifold $N$ of codimension
	$2$. If $M$ is connected, so is $M \setminus N$.
\end{lemma}
\begin{proof}
	Consider two points $x,y \in M \setminus N$. Since $M$ is connected, which
	is equivalent to be path-connected for manifolds, consider a path
	$\gamma : [0,1] \rightarrow M$ from $x$ to $y$. Then by theorem
	\cite[thm. 6.2.13]{Mukherjee-15}, $\gamma $ is homotopic to a path
	$\gamma '$ from $x$ and $y$ that is transverse to the submanifold
	$N$. Since as a real manifold, $\text{codim}(N) = 2$, the transversality
	with $\gamma '$ implies that
	$N \cap \text{im}(\gamma ') = \emptyset $, which shows that
	$M \setminus N$ is indeed path-connected.
\end{proof}

It is important to note that recent literature, like~\cite{Flatto-2009}
and~\cite{Dragovic-Radnovic-2025}, has explored Poncelet porisms in singular
cases, providing analytic conditions for configurations where conics are
tangent or degenerate. Our results complement these studies by characterizing
the global topology of the \textit{regular} Cayley set (non-degenerate conics
meeting transversely). From this perspective, the singular cases analyzed
in their work represent the singular locus that lies on the
\textbf{boundary} of the 9-dimensional manifold constructed herein, rather
than within the manifold itself.

\section{The moduli space of pairs of conics}
\label{sec::pair_of_conics}
	
In this section, we will construct explicitly a moduli space of pairs of
conics that intersect transversely. In the next section, we shall show
how this relates to elliptic curves. This relation is quite natural, since
as shown in~\cite{Dolgachev-12}, an elliptic curve defines a double cover
of $\mathbb C\mathbb P^{1}$, with four ramification points, which can be
considered as the base points of a pencil of conics. However, the construction
we propose here, does not seem to appear explicitly in the literature.
This construction is tailored to the context of Poncelet porism, through
the approach introduced in~\cite{griffiths_harris_78}, which relies on
elliptic curves.

Let us consider again pairs of non-degenerate conics that intersect transversely.
As in section~\ref{sec::manifold}, let $\mathfrak{O}$ be the open set of
$\mathbb C\mathbb P^{5} \times \mathbb C\mathbb P^{5}$ made of such pairs.

We define an action of the group of projective transformations
$G = \mathbb PGL_{3}(\mathbb C)$ of $\mathbb C\mathbb P^{2}$ on
$\mathfrak{O}$, as follows. Two pairs $(C,D)$ and $(C',D')$ are said to
be equivalent via congruence if we have: $C' \sim A^{T} C A$ and
$D' \sim A^{T} D A$ for some $A \in G$, where $\sim $ denotes equality
modulo multiplication by a non-vanishing scalar.

Notice that the pair
$(I_{3}, \text{diag}(\lambda _{1},\lambda _{2},\lambda _{3}))$ is equivalent via congruence to the pair
$(I_{3}, \text{diag}(\lambda _{\sigma (1)},\lambda _{\sigma (2)},\allowbreak
\lambda _{\sigma (3)}))$, for any permutation $\sigma \in S_{3}$. It appears
clearly by considering the permutation matrix
$A = (a_{ij}) \in GL_{3}(\mathbb C)$ defined by $\sigma $ (i.e.
$a_{ij} = \delta _{i,\sigma (j)}$). Therefore to each orbit, corresponds
a unique point in the quotient of projective plane by the third symmetric
group, $\mathbb C\mathbb P^{2} / S_{3}$. A more precise statement will
be presented below.

Now our first concern will be to understand what the isotropy group at
different points is.
	
\begin{proposition}
	%%LEAP%%%\label{prop2}
	\label{prop::isotropy_groups}
	Consider a pair $(C,D)$ equivalent to
	$(I_{3},\text{diag}(\lambda _{1},\lambda _{2},\lambda _{3}))$. If the point
	$[\lambda _{1}: \lambda _{2}: \lambda _{3}]$ in
	$\mathbb C\mathbb P^{2}/S_{3}$, is different from
	$[1:\omega :\omega ^{2}]$, where $\omega = e^{2\pi i/3}$, then the isotropy
	group is isomorphic to the Klein four-group $V_{4}$. At
	$(I_{3},\text{diag}(1,\omega ,\omega ^{2}))$, the isotropic group has 12
	elements and is isomorphic to $V_{4} \ltimes _{\Phi }\mathbb Z_{3}$, where
	$\Phi : \mathbb Z_{3} \rightarrow \text{Aut}(V_{4})$ is defined in the proof.
\end{proposition}
\begin{proof}
	It is well known that the isotropy groups are conjugate each other along
	each orbit. Therefore it is sufficient to compute these groups at points
	of the form
	$(I_{3},\text{diag}(\lambda _{1},\lambda _{2},\lambda _{3}))$. So let us
	consider a non-singular matrix $A$ such that
	$(A^{T} A, A^{T} \text{diag}(\lambda _{1},\lambda _{2},\lambda _{3}) A)$
	is equivalent to
	$(I_{3},\text{diag}(\lambda _{1},\lambda _{2},\lambda _{3}))$. Then there
	exists non-vanishing $\alpha $ and $\beta $, such that (i)
	$A^{T} A = \alpha I_{3}$ and (ii) $A^{T} \Lambda A = \beta \Lambda $, where
	$\Lambda = \text{diag}(\lambda _{1},\lambda _{2},\lambda _{3})$.
	
	From (i), we get $A^{T} = \alpha A^{-1}$, so that
	$\Lambda A = \frac{\beta}{\alpha} A \Lambda $. The standard basis
	$(e_{1},e_{2},e_{3})$ of $\mathbb C^{3}$ forms a basis of eigenvectors of
	$\Lambda $. Then we get:
	$\Lambda A e_{i} = \frac{\beta}{\alpha} \lambda _{i} A e_{i}$. Therefore
	there exists a permutation $\sigma \in S_{3}$ such that we have:
	$\frac{\beta}{\alpha} \lambda _{i} = \lambda _{\sigma (i)}$ and
	$A e_{i} = \mu _{i} e_{\sigma (i)}$ for all $i$ and for some
	$(\mu _{1},\mu _{2},\mu _{3}) \in (\mathbb C^{*})^{3}$. If $\sigma $ has
	a fixed point, then $\beta = \alpha $ and $\sigma = Id$ since the eigenvalues
	$\lambda _{i}$ are all distinct. Otherwise $\sigma $ is a cycle, so that
	$\sigma ^{3} = Id$. Therefore
	$\left (\frac{\beta}{\alpha} \right )^{3} = 1$ and
	$\frac{\beta}{\alpha} \in \{1, \omega , \omega ^{2}\}$.
	
	All together, if $\frac{\beta}{\alpha} = 1$, the matrix $A$ is diagonal:
	$A = \text{diag}(\mu _{1}, \mu _{2}, \mu _{3})$. In that case,
	$\mu _{i}^{2} = 1$, so that $A$ is either
	$I_{3},\text{diag}(1,1,-1) \sim \text{diag}(-1,-1,1)$ or
	$\text{diag}(1,-1,1) \sim \text{diag}(-1,1,-1)$, $\text{diag}(1,-1,-1)$. This shows that the isotropy group is indeed isomorphic
	to the Klein four-group $V_{4}$ (see~\cite[p.20]{Barnes-63}).
	
	If $\frac{\beta}{\alpha}$ is either $\omega $ or $\omega ^{2}$, then
	$[\lambda _{1}: \lambda _{2}: \lambda _{3}]$ is
	$[1:\omega :\omega ^{2}]$ up to a permutation. In that case
	$A = (a_{ij})$ is a kind of permutation matrix, where
	$a_{ij} = \mu _{j} \delta _{i,\sigma (j)}$. Then we have
	$\mu _{j}^{2} = 1$ for all $j$. Therefore the isotropy group has 12 elements,
	which are $P_{\sigma }A$ for $A \in V_{4}$ and $P_{\sigma}$ is the permutation
	matrix associated to the one of the 3-cycles, $Id,(123),(132)$. These cycles
	form a subgroup of $S_{3}$, isomorphic to $\mathbb Z_{3}$.
	
	By unfolding the product, we find that
	$P_{\sigma _{1}} \Lambda P_{\sigma _{2}} \Lambda ' = P_{\sigma _{1}
		\sigma _{2}} \Phi (\sigma _{2})(\Lambda ) \Lambda '$, where
	$\Phi : \mathbb Z_{3} \rightarrow \text{Aut}(V_{4}), \sigma \mapsto (
	\text{diag}(\lambda _{1},\lambda _{2},\lambda _{3}) \mapsto
	\text{diag}(\lambda _{\sigma (1)},\lambda _{\sigma (2)},\lambda _{
		\sigma (3)}))$. This defines an isomorphism from the group
	$V_{4} \ltimes _{\Phi }\mathbb Z_{3}$ to the isotropy group, where the
	product in $V_{4} \ltimes _{\Phi }\mathbb Z_{3}$ is defined by
	$(\Lambda ,\sigma _{1}).(\Lambda ',\sigma _{2}) = (\Phi (\sigma _{2})(
	\Lambda )\Lambda ',\sigma _{1} \sigma _{2})$.
\end{proof}
	
\textit{Remark:} A more geometric argument shows that if $A$ is in the isotropy
group of a pair $(C,D)$, then it is a kind of permutation matrix the elements
of which are all either $-1$ or $1$. Indeed, $C$ and $D$ are non-degenerate
conics that intersect in four distinct points:
$P = \{p_{1},p_{2},p_{3},p_{4}\}$. Since neither $C$ nor $D$ is degenerate,
the points in $P$ are in general position. Let $A \in G$, such that
$A^{T} C A \sim C$ and $A^{T} D A \sim D$, where $\sim $ means as above
equality modulo multiplication by a non-zero scalar. Then the set
$P$ is globally invariant by $A$ and every conics in the pencil defined
by $C$ and $D$ is globally invariant by $A$. Therefore the conics defined
by pairs of lines $(\overline{p_{1}p_{2}},\overline{p_{3}p_{4}})$,
$(\overline{p_{1}p_{3}},\overline{p_{2}p_{4}})$ and
$(\overline{p_{1}p_{4}},\overline{p_{2}p_{3}})$ are each one globally invariant
by $A$. Therefore the points
$q_{1} = \overline{p_{1}p_{2}} \cap \overline{p_{3}p_{4}}, q_{2} =
\overline{p_{1}p_{3}} \cap \overline{p_{2}p_{4}}, q_{3} =
\overline{p_{1}p_{4}} \cap \overline{p_{2}p_{3}}$ are fixed points of
$A$ and are in general position. Up to some 3-cycle, the reference frame
can be chosen so that $q_{1} =[1:0:0], q_{2}=[0:1:0], q_{3}=[0:0:1]$. In
that frame, $A$ is a diagonal matrix
$A = \text{diag}(\alpha _{1},\alpha _{2},\alpha _{3})$, with
$\alpha _{i}^{2} = 1$ for each $i$. Finally, since the numbering of the
points is arbitrary, $A$ is of the form
$P_{\sigma }\text{diag}(\alpha _{1},\alpha _{2},\alpha _{3})$, for
$\sigma \in \{Id,(123),(132)\}$. To get a finer result, one would use the
more algebraic approach we present in the proof of Proposition~\ref{prop::isotropy_groups}.

Now we shall prove that the orbits are embedded submanifolds of
$\mathfrak{O}$ of dimension $8$.

\begin{theorem}
	\label{thm::orbits}
	Each orbit is an embedded connected complex submanifold of
	$\mathfrak{O}$ of dimension $8$. More precisely the orbit of
	$(I_{3},\text{diag}(\lambda _{1},\lambda _{2},\lambda _{3}))$, where
	$[\lambda _{1}: \lambda _{2} : \lambda _{3}] \neq [1:\omega :\omega ^{2}]$
	in $\mathbb C\mathbb P^{2} /S_{3}$ is biholomorphic to $G/V_{4}$, while
	the orbit of $(I_{3},\text{diag}(1,\omega ,\omega ^{2}))$ is biholomorphic
	to $G/(V_{4} \ltimes _{\Phi }\mathbb Z_{3})$.
	
	We shall call the orbit of
	$(I_{3},\text{diag}(1,\omega ,\omega ^{2}))$
	\textnormal{the special orbit}, while the other orbits will be called
	\textnormal{standard orbits}.
\end{theorem}
\begin{proof}
	Let $(C,D) \in \mathfrak{O}$. Assume it is equivalent to
	$(I_{3},\text{diag}(\lambda _{1},\lambda _{2},\lambda _{3}))$, where
	$[\lambda _{1}: \lambda _{2} : \lambda _{3}] \neq [1:\omega :\omega ^{2}]$
	in $\mathbb C\mathbb P^{2} /S_{3}$. Then the isotropy group of
	$(C,D)$ is $V_{4}$ according to Proposition~\ref{prop::isotropy_groups}.
	
	By theorem~\cite[th.3.58]{Warner-83}, the quotient space $G/V_{4}$ is a
	complex manifold of dimension
	$\dim (G) - \dim (V_{4}) = \dim (G) = 8$. Then by theorem~\cite[th.3.62]{Warner-83},
	the orbit of $(C,D)$ is biholomorphic to $G/V_{4}$ and therefore has dimension
	$8$.
	
	The same argument is valid for the orbit of
	$(I_{3},\text{diag}(1,\omega ,\omega ^{2}))$, since the isotropy group,
	$V_{4} \ltimes _{\Phi }\mathbb Z_{3}$ is finite in this case too.
\end{proof}
	
In order to grasp more concretely the difference between the standard orbits
and the special one, we will make explicit a major difference between their
respective fundamental groups.

%p3 #&#
\begin{proposition}
	%%LEAP%%%\label{prop3}
	\label{prop::fund_group_orbits}
	The fundamental group of $G = \mathbb PGL_{3}(\mathbb C)$ is
	$\mathbb Z_{3}$, so that the fundamental group of the standard orbits has
	order $12$ and is isomorphic to $V_{4} \times \mathbb Z_{3}$, while for
	the special orbit, it has order $36$.
\end{proposition}
\begin{proof}
	The fundamental group of $G = \mathbb PGL_{3}(\mathbb C)$ results from
	a standard computation, that we recall schematically.
	
	Let
	$\eta : GL_{3}(\mathbb C) \rightarrow \mathbb PGL_{3}(\mathbb C)$ be the
	canonical projection. We have the following exact short sequence of Lie
	groups:
	\begin{equation*}
		1 \rightarrow \mathbb C^{*} \overset{i}{\rightarrow} GL_{3}(\mathbb C)
		\overset{\eta}{\rightarrow} \mathbb PGL_{3}(\mathbb C) \rightarrow 1,
	\end{equation*}
	where $i(\lambda ) = \lambda I_{3}$. This shows that $\eta $ is actually
	a fiber bundle, where each fiber is isomorphic to $\mathbb C^{*}$. Therefore
	by theorems~\cite[th. 4.41, th. 4.48]{hatcher-02}, we have a long exact
	sequence of homotopy groups:
	%
	%e4 #&#
	\begin{equation}
		\label{long_seq_PGL}
		\begin{aligned}
			\ldots \rightarrow{}& \pi _{2}(\mathbb C^{*}) \rightarrow \pi _{2}(GL_{3}(
			\mathbb C)) \rightarrow \pi _{2}(\mathbb PGL_{3}(\mathbb C))
			\rightarrow
			\\
			&\pi _{1}(\mathbb C^{*}) \overset{i_{*}}{\rightarrow} \pi _{1}(GL_{3}(
			\mathbb C)) \rightarrow \pi _{1}(\mathbb PGL_{3}(\mathbb C))
			\rightarrow \pi _{0}(\mathbb C^{*}) \rightarrow \ldots
		\end{aligned}
		%
		%%LEAP%%%\label{eq4}
	\end{equation}
	
	Since $\mathbb C^{*}$ is path-connected,
	$\pi _{0}(\mathbb C^{*}) = \{0\}$. Moreover $\mathbb C^{*}$ is homotopy
	equivalent to $S^{1}$, so that
	$\pi _{1}(\mathbb C^{*}) = \mathbb Z$ and
	$\pi _{2}(\mathbb C^{*}) = \{0\}$.
	
	Therefore we need to compute the homotopy groups of
	$GL_{3}(\mathbb C)$ in order to deduce from that the fundamental group
	of $\mathbb PGL_{3}(\mathbb C)$.
	
	According to~\cite{Borel-52}, every connected Lie group is homotopy equivalent
	to a maximal compact subgroup. Applied to our case
	$GL_{3}(\mathbb C)$, we get that $GL_{3}(\mathbb C) \simeq U(3)$, the unitary
	group. Therefore $\pi _{k}(GL_{3}(\mathbb C)) \cong \pi _{k}(U(3))$ for
	$k \geq 1$.
	
	Hence, we will compute the homotopy groups of $U(3)$. Consider the determinant
	map $\det : U(3) \rightarrow S^{1}$, which gives rise to the following
	exact sequence:
	\begin{equation*}
		1 \rightarrow SU(3) \rightarrow U(3) \overset{\det}{\rightarrow} S^{1}
		\rightarrow 1.
	\end{equation*}
	Similarly as above, this yields a long exact sequence of homotopy groups:
	\begin{equation*}
		\ldots \rightarrow \pi _{2}(S^{1}) \rightarrow \pi _{1}(SU(3))
		\rightarrow \pi _{1}(U(3)) \overset{\det _{*}}{\rightarrow} \pi _{1}(S^{1})
		\rightarrow \pi _{0}(SU(3)) \rightarrow \ldots
	\end{equation*}
	
	Since $SU(3)$ is connected and simply connected (see~\cite{Hall-2015}),
	we have:
	\begin{equation*}
		1 \rightarrow \pi _{1}(U(3)) \overset{\det _{*}}{\rightarrow} \pi _{1}(S^{1})
		\cong \mathbb Z\rightarrow 1,
	\end{equation*}
	so that
	$\pi _{1}(GL_{3}(\mathbb C)) \cong \pi _{1}(U(3)) \cong \mathbb Z$.
	
	Therefore the sequence~\eqref{long_seq_PGL} becomes:
	\begin{equation*}
		1 \rightarrow \pi _{2}(\mathbb PGL_{3}(\mathbb C)) \rightarrow
		\mathbb Z\overset{i_{*}}{\rightarrow} \mathbb Z\rightarrow \pi _{1}(
		\mathbb PGL_{3}(\mathbb C)) \rightarrow 1.
	\end{equation*}
	
	This yields:
	$\pi _{1}(\mathbb PGL_{3}(\mathbb C)) \cong \text{coker}(i_{*}) =
	\mathbb Z/ \text{im}(i_{*})$ and
	$\pi _{2}(\mathbb PGL_{3}(\mathbb C)) \cong \ker (i_{*})$, where
	$i_{*}: \pi _{1}(\mathbb C^{*}) \cong \mathbb Z\rightarrow \pi _{1}(GL_{3}(
	\mathbb C)) \cong \mathbb Z$ is the map induced by
	$i: \lambda \mapsto \lambda I_{3}$.
	
	Consider the loop $\gamma (t) = e^{2\pi i t}$ in $\mathbb C^{*}$, which
	correspond to $1$ in $\pi _{1}(\mathbb C^{*}) = \mathbb Z$. Then
	$i \circ \gamma (t) = e^{2 \pi i t} I_{3} \in SU(3)$. Thus, we have
	$\det (i \circ \gamma (t)) = e^{6 \pi i t} \in S^{1}$, which corresponds
	to $3$ in $\pi _{1}(S^{1}) = \mathbb Z$. Hence
	$i_{*}: \mathbb Z\rightarrow \mathbb Z$ is defined by $i_{*}(1) = 3$. This
	yields $\ker (i_{*}) = \{0\}$ and $\text{im}(i_{*}) = 3\mathbb Z$. As a
	conclusion we get: $\pi _{2}(\mathbb PGL_{3}(\mathbb C)) = \{0\}$ and
	\begin{equation*}
		\pi _{1}(\mathbb PGL_{3}(\mathbb C)) = \mathbb Z/ (3\mathbb Z) =
		\mathbb Z_{3},
	\end{equation*}
	as expected.\goodbreak
	
	As for the fundamental groups of $G/V_{4}$ and
	$G/(V_{4} \ltimes_{\Phi} \mathbb Z_{3})$, notice that the two maps
	$G \rightarrow G/V_{4}$ and
	$G \rightarrow G/(V_{4} \ltimes_{\Phi} \mathbb Z_{3})$ are both covering maps of
	degrees respectively $4$ and $12$. By proposition~\cite[prop. 1.32]{hatcher-02},
	this implies that $\pi _{1}(G/V4)$ and
	$\pi _{1}(G/(V_{4} \ltimes _{\Phi }\mathbb Z_{3}))$ are extension of
	$\pi _{1}(G) = \mathbb Z_{3}$ and have respectively order $12$ and
	$36$.
	
	Now, we shall determine the exact form of $\pi _{1}(G/V4)$. The short sequence:
	$1 \rightarrow V_{4} \rightarrow G \overset{\eta}{\rightarrow} G/V_{4}
	\rightarrow 1$ gives rise, as above, to a long sequence of homotopy groups:
	\begin{equation*}
		\ldots \pi _{1}(V_{4}) \rightarrow \pi _{1}(G) \rightarrow \pi _{1}(G/V_{4})
		\rightarrow \pi _{0}(V_{4}) \rightarrow \pi _{0}(G) \rightarrow
		\ldots
	\end{equation*}
	
	Since, $G$ is connected and $V_4$ is finite, we have:
	$\pi _{1}(V_{4}) = \{0\}$, $\pi _{0}(V_{4}) = V_{4}$ and
	$\pi _{0}(G) = \{0\}$. Notice that since $V_{4}$ is a finite group,
	$\pi _{0}(V_{4})$ can be seen as a group isomorphic to $V_{4}$ itself.
	In this context, the map
	$\pi _{1}(G/V_{4}) \rightarrow \pi _{0}(V_{4}) = V_{4}$ is a group homomorphism.
	Therefore this yields:
	\begin{equation*}
		0 \rightarrow \mathbb Z_{3} \overset{\eta _{*}}{\rightarrow} \pi _{1}(G/V_{4})
		\rightarrow V_{4} \rightarrow 1
	\end{equation*}
	
	This shows that $\eta _{*}(\pi _{1}(G))$ is a normal subgroup of
	$\pi _{1}(G/V_{4})$ as expected since $V_{4}$ is finite.
	
	According to~\cite{Brown-82}, to prove that $\pi _{1}(G/V_{4}))$ is the
	direct product $\mathbb Z_{3} \times V_{4}$, we need to prove that (i)
	the action of $V_{4}$ by conjugation on $\eta _{*}(\mathbb Z_{3})$ is trivial
	and that (ii) The classifying cohomology class
	$[\alpha ] \in H^{2}(V_{4},Z(\pi _{1}(G)))$ is zero, where
	$Z(\pi _{1}(G))$ is the center of $\pi _{1}(G)$.
	
	To verify (i), recall that the action of $V_4$ on $\eta_*(\pi_1(G))$ appearing in
	this extension is realized by the deck transformations of the covering
	$\eta : G \to G/V_4$, that is, by the translations $R_h : g \mapsto gh$ for
	$h \in V_4$ (see~\cite[\S1.3]{hatcher-02}). Since $G = \mathbb{P}GL_3(\mathbb{C})$
	is connected, any $h \in V_4$ can be joined to the identity by a path
	$c : [0,1] \to G$ with $c(0) = I_3$ and $c(1) = h$. The map
	$H : G \times [0,1] \to G$, $H(g,s) = g\,c(s)$, is then a homotopy from
	$\mathrm{id}_G = R_{I_3}$ to $R_h$; hence $(R_h)_* = \mathrm{id}$ on $\pi_1(G)$,
	so each $h$ acts trivially. As $\pi_1(G) = \mathbb{Z}_3$ is abelian, the fundamental group $\pi_1(G)$ the basepoint change ambiguity along $c$ (an inner automorphism of $\pi_1(G)$) is trivial as
	well. Therefore the conjugation action $\rho : V_4 \to \mathrm{Aut}(\pi_1(G))$ is
	trivial.

	As for the second cohomology group $H^{2}(V_{4},\mathbb Z_{3})$, it is
	trivial since the orders of $\mathbb Z_{3}$ and $V_{4}$ are relatively
	prime. Indeed the group $H^{2}(V_{4},\mathbb Z_{3})$ is annihilated by
	both~\cite[Cor. 10.2, p.84]{Brown-82} $|V_{4}|$ and
	$|\mathbb Z_{3}|$ and since they are relatively prime, the group is trivial.
	See more explicitly~\cite[Cor. 37,38, p.828]{Dummit-Foote-2004}
	
	Therefore the fundamental group $\pi _{1}(G/V_{4})$ is a trivial extension
	of $\pi _{1}(G)$ and we have:
	$\pi _{1}(G/V_{4}) = V_{4} \times \mathbb Z_{3}$.
\end{proof}
	
Now let us focus on the space of orbits
$\mathcal{M} = \mathfrak{O} / G$. Let
$\pi : \mathfrak{O} \rightarrow \mathcal{M}$ be the canonical projection.
As we have mentioned above, there is an injection from $\mathcal{M}$ into
the quotient $\mathbb C\mathbb P^{2}/S_{3}$. Now, we shall provide a more
precise statement, as follows.

%p4 #&#
\begin{proposition}
%%LEAP%%%\label{prop4}
\label{prop::moduli}
The space of orbits $\mathcal{M}$ can be endowed with the structure of
a two-dimensional complex manifold, such that it is isomorphic to
\begin{enumerate}
	\item an open set of the quotient $\mathbb C\mathbb P^{2}/S_{3}$, where
	$S_{3}$ is the third symmetric group,
	\item an open set of the weighted projective space
	$\mathbb P(1,2,3)$.
\end{enumerate}
In both cases, the open set is made of points which coordinates are distinct
and non-vanishing. In particular, $\mathcal{M}$ is connected.

Moreover the canonical projection $\pi $ is a holomorphic surjective submersion.
\end{proposition}
\begin{proof}
Let $(C,D) \in \mathfrak{O}$. Since these two conics are non-degenerate
and intersect transversely, the polynomial $\det (\lambda C + D)$ have
three distinct roots $\lambda _{1}, \lambda _{2}, \lambda _{3}$ and the
pair $(C,D)$ is equivalent via congruence to
$(I_{3},\text{diag}(\lambda _{1}, \lambda _{2}, \lambda _{3}))$.
\begin{enumerate}
	\item Let $\tilde{U} \subset \mathbb C\mathbb P^{2}$ be the open set defined
	by $\tilde{U} = \cap _{i=1}^{3} U_{i}$, where
	$U_{i} = \{[x_{0}:x_{1}:x_{2}] \in \mathbb C\mathbb P^{2}, x_{i}
	\neq 0\}$. Let $U$ be the open subset of $\tilde{U}$ made of points which
	coordinates of all distinct.
	
	Therefore each orbit or equivalently each element of the quotient space
	$\mathcal{M}$ is mapped via $\pi $ to the orbit of a point
	$[\lambda _{1} : \lambda _{2} : \lambda _{3}] \in U \subset \mathbb C
	\mathbb P^{2}$ under the action of the symmetric group $S_{3}$, since any
	permutation of the roots will define the same orbit in $\mathfrak{O}$.
	
	The discrete group $S_{3}$ acts on $U$ by permutation of the coordinates.
	Therefore it acts on $U$ holomorphically. The action is free since for
	$\sigma \in S_{3}$,
	$\sigma .[\lambda _{1} : \lambda _{2} : \lambda _{3}] = [\lambda _{
		\sigma (1)} : \lambda _{\sigma (2)} : \lambda _{\sigma (3)}]$ equals
	$[\lambda _{1} : \lambda _{2} : \lambda _{3}]$ only if $\sigma $ is the
	identity, since the $\lambda _{i}$ are distinct non-zero complex numbers.
	
	The action is also proper, i.e. the map
	$\phi : S_{3} \times U \rightarrow U \times U, (\sigma ,p) \mapsto (
	\sigma .p,p)$ is proper. Indeed consider a compact set
	$K \in U \times U$. Let $K_{1}$ and $K_{2}$ be respectively the projection
	of $K$ on the first factor and on the second factor. Then $K_{1}$ and
	$K_{2}$ are compact sets and $K \subset K_{1} \times K_{2}$. Then
	$\phi ^{-1}(K_{1} \times K_{2}) = \{(\sigma ,p) \in S_{3} \times U | p
	\in K_{2} \cap \sigma ^{-1}.K_{1}\}$ is a compact set since $S_{3}$ is
	finite, and $\phi ^{-1}(K) \subset \phi ^{-1}(K_{1} \times K_{2})$ is also
	compact, being a closed subset of a Hausdorff compact set.
	
	Therefore by the Holomorphic Quotient Manifold theorem~\cite{lee-2024},
	$U/S_{3}$ is a complex manifold. It is two-dimensional, since the covering
	map $U \rightarrow U/S_{3}$ is finite. Moreover $U$ is connected by Lemma~\ref{lem::connectivity_codim2}
	because it is the complement in $\tilde{U}$ of
	$Z((\lambda _{1}-\lambda _{2})(\lambda _{1}-\lambda _{3})(\lambda _{2}-
	\lambda _{3}))$. Therefore $U/S_{3}$ is connected too.
	
	The set $\mathcal{M}$ is in bijection with $U/S_{3}$ and as so, can be
	endowed with the structure of two-dimensional complex manifolds, biholomorphic
	to $U/S_{3}$.
	\item The space $\mathbb C\mathbb P^{2}/S_{3}$ is isomorphic to the weighted
	projective space $\mathbb P(1,2,3)$, as partly explained in~\cite{harris-92}.
	This can be seen as follows. Consider the map
	$\phi : \mathbb C\mathbb P^{2}/S_{3} \rightarrow \mathbb P(1,2,3), [
	\lambda _{1}:\lambda _{2}:\lambda _{3}] \mapsto [\lambda _{1}+
	\lambda _{2}+\lambda _{3}:\lambda _{1} \lambda _{2} + \lambda _{1}
	\lambda _{3} + \lambda _{2} \lambda _{3}: \lambda _{1} \lambda _{2}
	\lambda _{3}] = [e_{1}:e_{2}:e_{3}]$. It is obviously well-defined and
	bijective. The inverse map is given by the fact that
	$\lambda _{1}, \lambda _{2}, \lambda _{3}$ are the roots of
	$t^{3} - e_{1} t^{2} + e_{2} t - e_{3} = 0$, and so are locally holomorphic
	functions of $e_{1},e_{2},e_{3}$.
	
	Thus the space $\mathcal{M}$ is isomorphic to the open set of
	$\mathbb P(1,2,3)$ made of points with non-zero distinct coordinates.
\end{enumerate}

Now, notice that the projection
$\pi : \mathfrak{O} \rightarrow \mathcal{M}$ is holomorphic. This can be
seen as follows. Since the roots $\lambda _{i}$ are distinct, the derivative
of $P(t) = \det (tC + D)$ does not vanish at any $\lambda _{i}$. By the
implicit function theorem, this implies that for each pair
$(C_{0},D_{0}) \in \mathfrak{O}$, there exists an open neighborhood
$V$ of $(C_{0},D_{0})$, such that each $\lambda _{i}$ is a holomorphic
function of the coefficients of $P(t)$ and as a consequence a holomorphic
function of the entries of $C$ and $D$ whenever $(C,D) \in V$. We shall
simply write
$\lambda _{i}:V \rightarrow \mathbb C, (C,D) \mapsto \lambda _{i}(C,D)$.
Therefore restricted to $V$, $\pi $ is given by either by
$\pi (C,D) = [\lambda _{1}:\lambda _{2}:\lambda _{3}] \in \mathbb C
\mathbb P^{2}/S_{3}$ or by
$\pi (C,D) = [e_{1}:e_{2}:e_{3}] \in \mathbb P(1,2,3)$, where
$e_{1},e_{2},e_{3}$ are the symmetric functions of
$\lambda _{1}(C,D),\lambda _{2}(C,D),\lambda _{3}(C,D)$, as introduced
just above.

Let us now prove that $\pi $ is a submersion. This follows from the expression
of the derivatives of $\lambda _{i}$. We shall prove the following lemma.
		
%l3 #&#
\begin{lemma}
	\label{lem3}
	Let $(C, D)$ be a pair of non-degenerate conics that intersect transversely.
	Let $\lambda $ be a simple root of $P(x) = \det (xC +D)$. Let $v$ be a
	corresponding eigenvector satisfying $(\lambda C + D)v = 0$. Then the gradient
	of $\lambda $ with respect to the entries of $C$ and $D$ is given by the
	symmetric matrices:
	\begin{equation*}
		\nabla _{D} \lambda = - \frac{1}{v^{T} C v} vv^{T} \quad \text{and}
		\quad \nabla _{C} \lambda = -\frac{\lambda}{v^{T} C v} vv^{T},
	\end{equation*}
	where $v$ generates the kernel of $\lambda C +D$, or equivalently
	$v$ is a eigenvector of $-C^{-1}D$ with respect to $\lambda $.
	
	Furthermore, the directional derivative of $\lambda $ along a path with
	tangent vector $(\dot{C}, \dot{D})$ is:
	\begin{equation*}
		\dot{\lambda} = -
		\frac{v^{T} \dot{D} v + \lambda (v^{T} \dot{C} v)}{v^{T} C v}
	\end{equation*}
\end{lemma}
	
\begin{proof}
	The eigenvalue $\lambda $ is defined implicitly by the equation
	$P(\lambda ) = \det (\lambda C + D) = 0$. Since $\lambda $ is a simple
	root, $\frac{\partial P}{\partial \lambda} \neq 0$, so by the Implicit
	Function Theorem, we can compute the gradient of $\lambda $ with respect
	to a matrix entry $X_{ij}$ as
	$\frac{\partial \lambda}{\partial X_{ij}} = - (
	\frac{\partial P}{\partial X_{ij}}) / (
	\frac{\partial P}{\partial \lambda})$.
	
	We use Jacobi's formula (see~\cite{Lafontaine-2010}),
	$\frac{d}{dt}\det (A) = \text{tr}(\text{adj}(A)\frac{dA}{dt})$, to find the
	partial derivatives of $P$. Here $\text{adj}(A)$ is the adjoint matrix of
	$A$, i.e. the transpose of the cofactor matrix.
	\begin{enumerate}
		\item The derivative with respect to $\lambda $:
		\begin{equation*}
			\frac{\partial P}{\partial \lambda} = \text{tr}\left ( \text{adj}(
			\lambda C + D) \frac{\partial}{\partial \lambda}(\lambda C + D)
			\right ) = \text{tr}(\text{adj}(\lambda C + D)C)
		\end{equation*}
		\item The derivative with respect to an entry $D_{ij}$:
		\begin{equation*}
			\frac{\partial P}{\partial D_{ij}} = (\text{adj}(\lambda C + D))_{ji}
		\end{equation*}
		\item The derivative with respect to an entry $C_{ij}$:
		\begin{equation*}
			\frac{\partial P}{\partial C_{ij}} = \lambda (\text{adj}(\lambda C + D))_{ji}
		\end{equation*}
	\end{enumerate}
	Since $\lambda $ is a simple root, the matrix $\lambda C + D$ has rank
	2. Its adjoint is a rank-1 matrix proportional to $vv^{T}$, where
	$v$ is the eigenvector defined above such that $(\lambda C + D)v=0$. Let
	$\text{adj}(\lambda C + D) = k \cdot vv^{T}$ for some non-zero constant
	$k$. The denominator becomes:
	\begin{equation*}
		\frac{\partial P}{\partial \lambda} = \text{tr}(k \cdot vv^{T} C) = k
		\cdot \text{tr}(v^{T} C v) = k(v^{T} C v)
	\end{equation*}
	The term $v^{T} C v \neq 0$, as its vanishing would imply a non-transverse
	intersection of the conics. Indeed if $v^{T} C v = 0$, then also
	$v^{T} D v = 0$, so that $v$ is an intersection point and the tangent lines
	$Cv$ and $D v$ are then identical since $\lambda Cv = -Dv$, contradicting
	the transversality of the intersection.
	
	Now we assemble the gradients. The $(i,j)$-th entry of the gradient matrix
	$\nabla _{D} \lambda $ is:
	\begin{equation*}
		(\nabla _{D} \lambda )_{ij} =
		\frac{\partial \lambda}{\partial D_{ij}} = -
		\frac{k \cdot v_{j} v_{i}}{k(v^{T} C v)} = -
		\frac{v_{i} v_{j}}{v^{T} C v}
	\end{equation*}
	And for $\nabla _{C} \lambda $:
	\begin{equation*}
		(\nabla _{C} \lambda )_{ij} =
		\frac{\partial \lambda}{\partial C_{ij}} = -
		\frac{\lambda (k \cdot v_{j} v_{i})}{k(v^{T} C v)} = -\lambda
		\frac{v_{i} v_{j}}{v^{T} C v}
	\end{equation*}
	These expressions yield the matrices
	$\nabla _{D} \lambda = -\frac{1}{v^{T} C v} vv^{T}$ and
	$\nabla _{C} \lambda = -\frac{\lambda}{v^{T} C v} vv^{T}$.
	
	The directional derivative $\dot{\lambda}$ is the total derivative along
	the path $(\dot{C}, \dot{D})$:
	\begin{align*}
		\dot{\lambda} &= \text{tr}((\nabla _{C} \lambda )^{T} \dot{C}) +
		\text{tr}((\nabla _{D} \lambda )^{T} \dot{D})
		\\
		&= \text{tr}\left ( \left (-\frac{\lambda}{v^{T} C v} vv^{T}\right )
		\dot{C} \right ) + \text{tr}\left ( \left (-\frac{1}{v^{T} C v} vv^{T}
		\right ) \dot{D} \right )
		\\
		&= -\frac{1}{v^{T} C v} \left ( \lambda \cdot \text{tr}(vv^{T} \dot{C})
		+ \text{tr}(vv^{T} \dot{D}) \right )
	\end{align*}
	Since $\text{tr}(vv^{T} A) = v^{T} A v$, we arrive at the final expression:
	\begin{equation*}
		\dot{\lambda} = -
		\frac{v^{T} \dot{D} v + \lambda (v^{T} \dot{C} v)}{v^{T} C v}
	\end{equation*}
	This completes the proof.
\end{proof}
	
Now it remains to show that $\pi $ is a submersion. We will construct two
specific tangent vectors, $T_{1}$ and $T_{2}$, and show their images under
$d\pi $ are linearly independent, thus spanning the target tangent space.
To each root $\lambda _{i}$ is associated a eigenvector $v_{i}$, such that
$(\lambda _{i} C + D) v_{i} = 0$.

Let's define a tangent vector $T_{1}$ represented by the pair
$(\dot{C}_{1}, \dot{D}_{1}) = (0, C v_{1} v_{1}^{T} C)$. This is a valid
tangent vector as $\dot{C}_{1}$ and $\dot{D}_{1}$ are symmetric, and
$\dot{D}_{1}$ is not a multiple of $D$, since
$\text{rank}(D_{1}) = 1 \neq 3 = \text{rank}(D)$.

Let's compute the change in the eigenvalues induced by $T_{1}$. For
$\lambda _{1}$:
\begin{align*}
	\dot{\lambda}_{1} &= -
	\frac{v_{1}^{T} (C v_{1} v_{1}^{T} C) v_{1} + \lambda _{1} (v_{1}^{T} (0) v_{1})}{v_{1}^{T} C v_{1}}
	\\
	&= - \frac{(v_{1}^{T} C v_{1})(v_{1}^{T} C v_{1})}{v_{1}^{T} C v_{1}} =
	- v_{1}^{T} C v_{1}
\end{align*}
This is non-zero because transversality ensures
$v_{1}^{T} C v_{1} \neq 0$. For $\lambda _{2}$ and $\lambda _{3}$:
\begin{align*}
	\dot{\lambda}_{k} &= -
	\frac{v_{k}^{T} (C v_{1} v_{1}^{T} C) v_{k}}{v_{k}^{T} C v_{k}}
	\\
	&= - \frac{(v_{k}^{T} C v_{1})(v_{1}^{T} C v_{k})}{v_{k}^{T} C v_{k}}
	\quad (\text{for } k=2,3)
\end{align*}
Because the eigenvectors are $C$-orthogonal, $v_{k}^{T} C v_{1} = 0$ for
$k \neq 1$, since the roots $\lambda _{i}$ are distinct. Therefore,
$\dot{\lambda}_{2} = 0$ and $\dot{\lambda}_{3} = 0$. The change induced
by $T_{1}$ is
$(\dot{\lambda}_{1}, \dot{\lambda}_{2}, \dot{\lambda}_{3}) = (-v_{1}^{T}
C v_{1}, 0, 0)$.

Similarly, define $T_{2} = (0, -C v_{2} v_{2}^{T} C)$. An identical calculation
using the $C$-orthogonality shows that the change induced by $T_{2}$ is:
\begin{equation*}
	(\dot{\lambda}_{1}, \dot{\lambda}_{2}, \dot{\lambda}_{3}) = (0, v_{2}^{T}
	C v_{2}, 0)
\end{equation*}
where $v_{2}^{T} C v_{2} \neq 0$.

Now in the local coordinates
$(\lambda _{1}/\lambda _{3}, \lambda _{2}/\lambda _{3})$ of the codomain,
we get:
\begin{enumerate}
	\item for $T_{1}$, the change is
	\begin{equation*}
		\left (
		\frac{\dot{\lambda}_{1}\lambda _{3} - \lambda _{1}\dot{\lambda}_{3}}{\lambda _{3}^{2}},
		\frac{\dot{\lambda}_{2}\lambda _{3} - \lambda _{2}\dot{\lambda}_{3}}{\lambda _{3}^{2}}
		\right ) = \left (
		\frac{-(v_{1}^{T} C v_{1})\lambda _{3} - 0}{\lambda _{3}^{2}}, 0
		\right ) = \left (-\frac{v_{1}^{T} C v_{1}}{\lambda _{3}}, 0\right ).
	\end{equation*}
	\item for $T_{2}$, the change is
	\begin{equation*}
		\left (0,
		\frac{(v_{2}^{T} C v_{2})\lambda _{3} - 0}{\lambda _{3}^{2}}\right ) =
		\left (0, \frac{v_{2}^{T} C v_{2}}{\lambda _{3}}\right ).
	\end{equation*}
\end{enumerate}
	
Since $\lambda _{3}$, $v_{1}^{T} C v_{1}$, and $v_{2}^{T} C v_{2}$ are
all non-zero, the images of $T_{1}$ and $T_{2}$ are two linearly independent
vectors that form a basis for the 2-dimensional tangent space of the codomain.

This shows that $d\pi $ is surjective at every pair
$(C,D) \in \mathfrak{O}$, i.e. $\pi $ is a submersion.
\end{proof}
	
Note that this submersion cannot define a fiber bundle, because the fibers
are not all biholomorphic due to the special orbit, while the base
$\mathcal{M}$ is connected.

%s5 #&#
\section{Elliptic curves and pairs of conics}
%%LEAP%%%\label{sec5}
\label{sec::elliptic_curves}

As well known, an elliptic curve is a compact Riemann surface of genus
1 and every elliptic curve is isomorphic to an algebraic curve defined in
$\mathbb C^{2}$ by an equation of the form $y^{2} = P(x)$ where $P$ is
a polynomial of degree $3$ with three distinct roots:
$\lambda _{1}, \lambda _{2}, \lambda _{3}$. Through a change of coordinates,
we can always assume that each $\lambda _{i} \in \mathbb C^{*}$.

The relation between elliptic curves and linear pencil of curves has been
thoroughly investigated in the case of pencil of cubics. This is the so-called
Hesse pencil. See~\cite{Dolgachev-12} for details. While the relation to
pencil of conics is very natural too, it seems it has not been studied
in a comprehensive manner. We try to give here a more systematic treatment.
	
\begin{lemma}
	%%LEAP%%%\label{lem4}
	\label{lem::elliptic_curves_conics_pencil}
	\begin{enumerate}
		\item To every elliptic curve is associated a pair of conics $(C,D)$ in
		$\mathfrak{O}$, such that the curve is defined by
		$y^{2} = \det (xC+D)$.
		\item Every pair of conics $(C,D)$ in $\mathfrak{O}$ defines an elliptic
		curve by the cubic equation: $y^{2} = \det (x C + D)$.
	\end{enumerate}
\end{lemma}

The proof is standard, but for the sake of completeness, we shall present
here a possible version of it.
	
\begin{proof}
	\begin{enumerate}
		\item Consider the elliptic curve defined by $y^{2} = P(x)$, with three
		non-zero distinct roots $\lambda _{1}, \lambda _{2}, \lambda _{3}$. Then
		$P(x) = \det (xI_{3} - \text{diag}(\lambda _{1},\lambda _{2},\lambda _{3}))$.
		The pair of conics
		$(I_{3},- \text{diag}(\lambda _{1},\lambda _{2},\lambda _{3}))$ is indeed
		in $\mathfrak{O}$.
		\item The curve is smooth, since the roots of $\det (xC+D)$ are distinct.
		After homogenization, the equation becomes:
		$y^{2} z = \det (xC + zD)$. Let $E$ be the curve in
		$\mathbb C\mathbb P^{2}$ defined by this equation. Let the line at infinity
		be defined by $z=0$.
		
		In the affine piece defined by $z \neq 0$, we have a morphism
		$\xi : E_{z} = \{(x,y) \in \mathbb C^{2} | y^{2} = \det (xC+D)\}
		\rightarrow \mathbb C\mathbb P^{1}, (x,y) \mapsto [x:1]$. This defines
		a covering of degree $2$ with three ramification points of index 2 (that
		is $(\lambda _{1},0),(\lambda _{2},0),(\lambda _{3},0)$). The curve
		$E$ has a single point at infinity, which is $[0:1:0]$. We extend the map
		$\xi $, by sending this point at infinity to $[1:0]$. Then $[0:1:0]$ becomes
		also a ramification point of index $2$.
		
		The genus can then be computed through Hurwitz formula:
		$2g -2 = 2(2g(\mathbb C\mathbb P^{1})-2) + 4(2-1) = -4+4 = 0$, so that
		$g=1$ as expected.\qedhere
	\end{enumerate}
\end{proof}
	
The j-invariant of an elliptic curve $E$ defined by an equation
$y^{2} = x^{3} + ax^{2} + bx + c$ over a field $K$ is defined by (according
to~\cite{silverman-09} or \cite{Kunz-05}):
\begin{equation*}
	j(E) = -
	\frac{256 \left (a^{2} - 3 b\right )^{3}}{4 a^{3} c - a^{2} b^{2} - 18 a b c + 4 b^{3} + 27 c^{2}}
\end{equation*}

The denominator of $j(E)$ is the discriminant of
$x^{3} + ax^{2} + bx + c$ which does not vanish, since all roots are simple.

It is well known that the j-invariant classifies complex elliptic curves
up to isomorphism. See~\cite{silverman-09} or \cite{Kunz-05} for details.

%l5 #&#
\begin{lemma}
	\label{lem5}
	The j-invariant of an elliptic curve defined by a pair of conics
	$(C,D) \in \mathfrak{O}$ is given by:
	%
	%e5 #&#
	\begin{equation}
		\label{eq::jinv_omega}
		j(C,D) =
		\frac{256 (3 \sigma _{12} \sigma _{30} - \sigma _{21}^{2})^{3}}{\sigma _{30}^{2} (27 \sigma _{03}^{2} \sigma _{30}^{2} + 2 \sigma _{12} \sigma _{30} (- 9 \sigma _{03} \sigma _{21} + 2 \sigma _{12}^{2}) + \sigma _{21}^{2} (4 \sigma _{03} \sigma _{21} - \sigma _{12}^{2}))}
		%%LEAP%%%\label{eq5}
	\end{equation}
	It is therefore a holomorphic function on $\mathfrak{O}$.
\end{lemma}
\begin{proof}
	A cubic curve defined by
	$y^{2} = \alpha x^{3} + \beta x^{2} + \gamma x + \delta $ (with
	$\alpha \neq 0$ obviously) is isomorphic to the curve defined by
	$y^{2} = x^{3} + \frac{\beta}{\alpha} x^{2} + \frac{\gamma}{\alpha} x +
	\frac{\delta}{\alpha}$, by the isomorphism
	$(x,y) \mapsto (x,y/\sqrt{\alpha})$, for any square root of
	$\alpha $.
	
	Therefore if we consider a pair of conics $(C,D) \in \mathfrak{O}$. Then
	the curve
	$y^{2} = \det (xC+D) = \sigma _{3,0} x^{3} + \sigma _{2,1} x^{2} +
	\sigma _{1,2} x + \sigma _{0,3}$ is isomorphic to the curve defined by
	$y^{2} = x^{3} + \frac{\sigma _{2,1}}{\sigma _{3,0}} x^{2} +
	\frac{\sigma _{1,2}}{\sigma _{3,0}} x +
	\frac{\sigma _{0,3}}{\sigma _{3,0}}$. Therefore the j-invariant is given
	by the expression~\eqref{eq::jinv_omega}, which is obviously a holomorphic
	function defined on $\mathfrak{O}$.
\end{proof}
	
\begin{theorem}
	\label{thm4}
	The j-invariant on $\mathfrak{O}$ factors through the projection
	$\pi : \mathfrak{O} \rightarrow \mathcal{M}$:
	\begin{equation*}
		\begin{tikzcd}
			\mathfrak{O} \arrow[r, "j"] \arrow[d,"\pi "'] & \mathbb C
			\\
			\mathcal{M} \arrow[ru, "\overline{j}"']
		\end{tikzcd}
	\end{equation*}
		
	Moreover, $\overline{j}$ is a submersion outside the curve
	$X \subset \mathcal{M} \subset \mathbb C\mathbb P^{2}/S_{3}$ defined by
	the equation
	%
	%e6 #&#
	\begin{equation}
		\label{eq::critical_points_overline_j}
		\left (\lambda _{1} - 2 \lambda _{2} + \lambda _{3}\right ) \left (
		\lambda _{1} + \lambda _{2} - 2 \lambda _{3}\right ) \left (2
		\lambda _{1} - \lambda _{2} - \lambda _{3}\right ) \left (\lambda _{1}^{2}
		- \lambda _{1} \lambda _{2} - \lambda _{1} \lambda _{3} + \lambda _{2}^{2}
		- \lambda _{2} \lambda _{3} + \lambda _{3}^{2}\right )^{2} = 0
		%%LEAP%%%\label{eq6}
	\end{equation}
		
	Note that $\overline{j}$ is not related to complex conjugacy. The critical
	values of $\overline{j}$ and $j$ are $0,1728$, which correspond to curves
	having a special group of automorphisms.
	
	Finally the function $j$ is a submersion outside the inverse image
	$\pi ^{-1}(X)$.
\end{theorem}
\begin{proof}
	Two pairs of conics $(C,D)$ and $(C',D')$ in $\mathfrak{O}$ are equivalent
	if there exists a non-singular matrix $A$, such that
	$C' \sim A^{T} C A$ and $D' \sim A^{T} D A$. Then
	$\det (xC'+D') \sim \det (x A^{T} C A + A^{T} D A) = \det (A)^{2}
	\det (xC +D)$. Then the curves $y^{2} = \det (xC'+D')$ and
	$y^{2} = \det (xC+D)$ are isomorphic and as so, have the same j-invariant.
	
	Since $j$ is invariant over congruence equivalence classes, it factors
	through $\pi $. It is surjective since for every complex number, there
	are elliptic curves having this number as their j-invariant as shown in~\cite{Kunz-05,silverman-09}
	and by Lemma~\ref{lem::elliptic_curves_conics_pencil}, every elliptic curve
	can be defined by a pair of conics in $\mathfrak{O}$.
	
	As for $\overline{j}$, notice that it can be written explicitly as a function
	of $\lambda _{1}, \lambda _{2}, \lambda _{3}$:
	%
	%e7 #&#
	\begin{equation}
		\label{eq::j_lambda}
		\overline{j}(\lambda _{1},\lambda _{2},\lambda _{3}) =
		\frac{256 \left (\lambda _{1}^{2} - \lambda _{1} \lambda _{2} - \lambda _{1} \lambda _{3} + \lambda _{2}^{2} - \lambda _{2} \lambda _{3} + \lambda _{3}^{2}\right )^{3}}{\left (\lambda _{1} - \lambda _{2}\right )^{2} \left (\lambda _{1} - \lambda _{3}\right )^{2} \left (\lambda _{2} - \lambda _{3}\right )^{2}}
		%%LEAP%%%\label{eq7}
	\end{equation}
	
	This expression is merely obtained by substitution in equation~\eqref{eq::jinv_omega}.
	
	The gradient of $\overline{j}$ vanishes over the curve defined by
	$\left (\lambda _{1} - 2 \lambda _{2} + \lambda _{3}\right ) \left (
	\lambda _{1} + \lambda _{2} - 2 \lambda _{3}\right ) \left (2
	\lambda _{1} - \lambda _{2} - \lambda _{3}\right ) \left (\lambda _{1}^{2}
	-\right.\allowbreak \left. \lambda _{1} \lambda _{2} - \lambda _{1} \lambda _{3} + \lambda _{2}^{2}
	- \lambda _{2} \lambda _{3} + \lambda _{3}^{2}\right )^{2} = 0$.
	
	Then one can compute the critical values, which turn to be $0$ and
	$1728$. More precisely, over
	$\left (\lambda _{1} - 2 \lambda _{2} + \lambda _{3}\right ) \left (
	\lambda _{1} + \lambda _{2} - 2 \lambda _{3}\right ) \left (2
	\lambda _{1} - \lambda _{2} - \lambda _{3}\right )= 0$, the value of
	$\overline{j}$ is $1728$. And over
	$\left (\lambda _{1}^{2} - \lambda _{1} \lambda _{2} - \lambda _{1}
	\lambda _{3} + \right.\allowbreak \left.\lambda _{2}^{2} - \lambda _{2} \lambda _{3} +
	\lambda _{3}^{2}\right )^{2} = 0$,  $\overline{j}$ vanishes.
	
	Outside the fibers over these critical values $\overline{j}$ is submersion.
	
	Then since $j = \overline{j} \circ \pi $ and $\pi $ is a surjective submersion,
	the function $j$ happens to be a submersion outside $\pi ^{-1}(X)$.
\end{proof}

To conclude this section, let us describe geometrically the points where
$\overline{j}$ and $j$ fail to be submersions.

The critical variety $X$ from the above theorem is the union of $5$ lines,
defined in the dual projective plane by
$l_{1} = [-2:1:1], l_{2} = [1:-2:1], l_{3} = [1:1:-2], l_{4} = \left [
\frac{3}{2} + \frac{\sqrt{3} i}{2}: - \frac{3}{2} +
\frac{\sqrt{3} i}{2}: - \sqrt{3} i \right ], l_{5} = \left [
\frac{3}{2} - \frac{\sqrt{3} i}{2}: - \frac{3}{2} -
\frac{\sqrt{3} i}{2}: \sqrt{3} i \right ]$. The lines $l_{4},l_{5}$ appear
with multiplicity $2$, while $l_{1},l_{2},l_{3}$ have multiplicity
$1$ in the dual curve of $X$. By lifting these lines to
$\mathfrak{O}$, we get the following description:
\begin{enumerate}
	\item Pairs of conics, with a vanishing j-invariant are congruent to
	$(I_{3},D)$ where $D$ is a diagonal matrix whose diagonal elements define
	a point in $\mathbb C\mathbb P^{2}$ that lies in either $l_{4}$ or
	$l_{5}$. These pairs, regarded as elliptic curves, have an automorphism
	group equal to $\mathbb Z/6\mathbb Z$ (see \cite{Kunz-05} or
	\cite{silverman-09}). Also notice that the special orbit (see Theorem~\ref{thm::orbits})
	lies in the fiber $j^{-1}(0)$.
	\item Pairs of conics, with a j-invariant equal to $1728$ are congruent
	to $(I_{3},D)$ where $D$ is a diagonal matrix whose one element on the
	diagonal is the average of the other two elements. For these pairs, the
	group of automorphism is $\mathbb Z/4\mathbb Z$ (again see
	\cite{Kunz-05} or \cite{silverman-09})
\end{enumerate}
	
\section{Cayley set and J-invariant}
\label{sec::cayley_set_j_inv}

At this point, we shall understand the relation between
$\mathfrak{C}$ and $\mathcal{M}$. First, let us check that the action of
$G = \mathbb PGL_{3}(\mathbb C)$ is well-defined on $\mathfrak{C}$.

\begin{proposition}
	\label{prop::projection_Caycley_set}
	The orbit of a pair $(C,D) \in \mathfrak{C}$ is included in
	$\mathfrak{C}$. As a consequence the projection of $\mathfrak{C}$ on
	$\mathcal{M}$ is a Zariski closed set and is defined by the equation, where
	the left-hand side is an irreducible polynomial:
	%
	%e8 #&#
	\begin{equation}
		\label{eq::Cayley_M}
		- \lambda _{1}^{2} \lambda _{2}^{2} + 2 \lambda _{1}^{2} \lambda _{2}
		\lambda _{3} - \lambda _{1}^{2} \lambda _{3}^{2} + 2 \lambda _{1}
		\lambda _{2}^{2} \lambda _{3} + 2 \lambda _{1} \lambda _{2} \lambda _{3}^{2}
		- \lambda _{2}^{2} \lambda _{3}^{2} = 0
		%%LEAP%%%\label{eq8}
	\end{equation}
	Moreover, this equation defines a smooth algebraic set in
	$\mathcal{M}$.
\end{proposition}
\begin{proof}
	We have: $\det (tA^{T} C A + A^{T} D A) = \det (A)^{2} \det (tC+D)$. So
	the coefficients of the Taylor expansion of $\sqrt{\det (tC+D)}$ at
	$t=0$ are just multiplied by a non-zero scalar. Therefore
	$(A^{T} C A, A^{T} D A) \in \mathfrak{C}$ if, and only if
	$(C,D) \in \mathfrak{C}$. Consequently, if a pair $(C,D)$ is equivalent
	to $(I_{3},\text{diag}(\lambda _{1},\lambda _{2},\lambda _{3}))$, it lies
	in $\mathfrak{C}$ if, and only if
	$(I_{3},\text{diag}(\lambda _{1},\lambda _{2},\lambda _{3})) \in
	\mathfrak{C}$.
	
	By substitution into equation~\eqref{eq::cayley_equation}, we get equation~\eqref{eq::Cayley_M}.
	
	The left hand side, say
	$Q(\lambda _{1}, \lambda _{2}, \lambda _{3})$ of equation~\eqref{eq::Cayley_M}
	is necessarily an irreducible polynomial. Otherwise, different irreducible
	components would intersect, since $\mathbb C$ is algebraically closed.
	These points would have also lift to self-intersections in
	$\mathfrak{C}$, contradicting the smoothness of $\mathfrak{C}$.
	
	Moreover, if we consider the system
	$Q = \frac{\partial Q}{\partial \lambda _{1}} =
	\frac{\partial Q}{\partial \lambda _{2}} =
	\frac{\partial Q}{\partial \lambda _{3}} = 0$, we find that the only solutions
	in $\mathbb C\mathbb P^{2}$ are $[1:0:0]$, $[0:1:0]$ and $[0:0:1]$, the projections
	of which in $\mathbb C\mathbb P^{2}/S_{3}$ are outside of
	$\mathcal{M}$. As a consequence the algebraic set defined by this equation
	is smooth in $\mathcal{M}$.
\end{proof}

The points where the restriction of $\overline{j}$ to
$\pi (\mathfrak{C})$ is not a submersion are therefore of two types:
\begin{itemize}
	\item the intersection of the curves defined by equations~\eqref{eq::Cayley_M}
	and~\eqref{eq::critical_points_overline_j},
	\item and the points where the tangent of $\pi (\mathfrak{C})$ is in the
	kernel of the differential $d\overline{j}$.
\end{itemize}

The critical points of the first type are the elements of
$S = S_{0} \cup S_{1728}$, where we have:
\begin{equation*}
	S_{0} = \overline{j}_{|\pi (\mathfrak{C})}^{-1}(0) = \left \{ \left [
	\frac{ \left (-1 - \sqrt{3} i\right )}{2}: \left (- \frac{1}{2} +
	\frac{\sqrt{3} i}{2}\right ): 1\  \right ], \left [
	\frac{ \left (-1 + \sqrt{3} i\right )}{2}: \left (- \frac{1}{2} -
	\frac{\sqrt{3} i}{2}\right ): 1 \right ] \right \}
\end{equation*}
and if the points in $\mathbb{C} \mathbb{P}^2$ lying over $S_{1728} = \overline{j}_{|\pi (\mathfrak{C})}^{-1}(1728)$ are given by:
\begin{align*}
	%\begin{array}{l}
	& \left \lbrace \left [ \left (1 - \sqrt{-3 + 2 \sqrt{3}}\right ):
	\left (\sqrt{-3 + 2 \sqrt{3}} + 1\right ): 1 \right ], \right . \
	\
	\left [ \left (1 - \sqrt{- 2 \sqrt{3} - 3}\right ): \left (1 + \sqrt{-
		2 \sqrt{3} - 3}\right ): 1 \right ], \
	\\
	&\left [ \left (1 + \sqrt{- 2 \sqrt{3} - 3}\right ): \left (1 - \sqrt{-
		2 \sqrt{3} - 3}\right ) : 1 \right ], \
	\
	\left [ \left (\sqrt{-3 + 2 \sqrt{3}} + 1\right ): \left (1 - \sqrt{-3
		+ 2 \sqrt{3}}\right ): 1 \right ],
	\\
	&\left [ \left (1 + \sqrt{3} + \sqrt{3 + 2 \sqrt{3}}\right ): \left (
	\frac{\sqrt{3}}{2} + 1 + \sqrt{\frac{3}{4} + \frac{\sqrt{3}}{2}}
	\right ): 1 \right ],
	\\
	&\left [ \left (- \sqrt{3} + 1 - \sqrt{3 - 2 \sqrt{3}}\right ): \left (-
	\frac{\sqrt{3}}{2} + 1 - \sqrt{\frac{3}{4} - \frac{\sqrt{3}}{2}}
	\right ): 1 \right ],
	\\
	&\left [ \left (- \sqrt{3} + 1 + \sqrt{3 - 2 \sqrt{3}}\right ): \left (-
	\frac{\sqrt{3}}{2} + 1 + \sqrt{\frac{3}{4} - \frac{\sqrt{3}}{2}}
	\right ): 1 \right ],
	\\
	&\left [
	\frac{ \left (- \sqrt{3} + 2 - \sqrt{3 - 2 \sqrt{3}}\right )}{2}:
	\left (- \sqrt{3} + 1 - \sqrt{3 - 2 \sqrt{3}}\right ): 1 \right ],
	\\
	&\left [
	\frac{ \left (- \sqrt{3} + 2 + \sqrt{3 - 2 \sqrt{3}}\right )}{2}:
	\left (- \sqrt{3} + 1 + \sqrt{3 - 2 \sqrt{3}}\right ): 1 \right ],
	\\
	&\left [
	\frac{ \left (\sqrt{3} + 2 + \sqrt{3 + 2 \sqrt{3}}\right )}{2}:
	\left (1 + \sqrt{3} + \sqrt{3 + 2 \sqrt{3}}\right ): 1 \right ],
	\\
	&\left [ \left (- \sqrt{3 + 2 \sqrt{3}} + 1 + \sqrt{3}\right ): \left (-
	\sqrt{\frac{3}{4} + \frac{\sqrt{3}}{2}} + \frac{\sqrt{3}}{2} + 1
	\right ): 1 \right ], \
	\\
	&\left . \left [
	\frac{ \left (- \sqrt{3 + 2 \sqrt{3}} + \sqrt{3} + 2\right )}{2}:
	\left (- \sqrt{3 + 2 \sqrt{3}} + 1 + \sqrt{3}\right ): 1 \right ]
	\right \rbrace
	%%\end{array}
	%
\end{align*}

If we project these points into $\mathcal{M}$, we get the two following orbits as follows:
$$
S_{1728} = \overline{j}_{|\pi (\mathfrak{C})}^{-1}(1728) =
\left\{
\begin{aligned}
	& \left[\,1 + \tfrac{\sqrt{3}}{2} + \tfrac{1}{2}\sqrt{3 + 2\sqrt{3}}
	\ :\ 1 + \sqrt{3} + \sqrt{3 + 2\sqrt{3}} \ :\ 1\,\right], \\
	& \left[\,1 - \tfrac{\sqrt{3}}{2} + \tfrac{i}{2}\sqrt{2\sqrt{3} - 3}
	\ :\ 1 - \sqrt{3} + i\sqrt{2\sqrt{3} - 3} \ :\ 1\,\right]
\end{aligned}
\right\}.
$$

The critical points of the second type are those where the gradient of
$\overline{j}$ and the tangent (viewed as a point in the dual plane) to
the curve defined by~\eqref{eq::Cayley_M} are proportional, which yields
$3$ algebraic equations on $\lambda _{1},\lambda _{2},\lambda _{3}$. The
zero set of this system clearly contains the critical curve $X$ defined
by~\eqref{eq::critical_points_overline_j}. In addition it contains points
are not in $\mathcal{M}$. They are the points that are defined by
$[0:0:1],[0:1:0]$ in $\mathbb C\mathbb P^{2}$.

From this analysis, one concludes that
\textsc{$\overline{j}_{|\pi (\mathfrak{C}) \setminus S}: \pi (\mathfrak{C}) \setminus S \rightarrow \mathbb C\setminus \{0,1728\}$ is a submersion}.
More precisely, we have the following proposition.

\begin{proposition}
	\label{prop::covering}
	The restriction of $\overline{j}$ on
	$\pi (\mathfrak{C}) \setminus S$ defines a holomorphic covering space of
	degree $4$ over $\mathbb C\setminus \{0,1728\}$.
\end{proposition}

\begin{proof}
	Let $z \in \mathbb C\setminus \{0,1728\}$, the fiber of
	$\overline{j}$ over $z$ is an algebraic curve in $\mathcal{M}$ defined
	by:
	\begin{equation}
		\label{eq::fiber_overline_j}
		256 \left (\lambda _{1}^{2} - \lambda _{1} \lambda _{2} - \lambda _{1}
		\lambda _{3} + \lambda _{2}^{2} - \lambda _{2} \lambda _{3} +
		\lambda _{3}^{2}\right )^{3} = z \left (\lambda _{1} - \lambda _{2}
		\right )^{2} \left (\lambda _{1} - \lambda _{3}\right )^{2} \left (
		\lambda _{2} - \lambda _{3}\right )^{2}
	\end{equation}
	
	This curve lies in $\mathbb C\mathbb P^{2} / S_{3}$, since it is defined
	by a polynomial which is invariant to variable permutations. It is the
	image by the canonical projection of a complex algebraic curve embedded
	in $\mathbb C\mathbb P^{2}$. A short computation shows that the fiber is
	always smooth in $\mathcal{M}$ for $z \neq 0,1728$.
	
	By Bezout theorem, each fiber intersect $\pi (\mathfrak{C})$ which is defined
	by equation~\eqref{eq::Cayley_M}, in $24$ points counted with multiplicities in $\mathbb{C} \mathbb{P}^2$.
	The system formed by equations~\eqref{eq::Cayley_M} and~\eqref{eq::fiber_overline_j}
	shows that if one of the $\lambda _{i}$ vanishes, the two other also vanish.
	Furthermore, if two of the $\lambda _{i}$ are identical the system implies
	that they are all zero. This implies that all of the $24$ intersection
	points have distinct and non-zero coordinates. Therefore they project into $\mathcal{M}$. Since the orbit size modulo $S_3$ is $6$, this yields $4$ orbits or points in $\mathcal{M}$. 
	
	These points do not lie in $S$, then the function $\overline{j}$ is a submersion at each of these points.
	
	If the fiber over $z$ and $\pi (\mathfrak{C})$ intersects transversely
	at some point
	$\lambda = [\lambda _{1}: \lambda _{2} : \lambda _{3}]$, then the tangents
	at $\lambda $ generate the tangent space of $\mathcal{M}$ at
	$\lambda $. In that case, since the differential of $\overline{j}$ along
	the fiber $\overline{j}^{-1}(z)$ vanishes, the tangent to
	$\pi (\mathfrak{C})$ is mapped surjectively onto
	$T_{z}(\mathbb C) \cong \mathbb C$ by $\overline{j}$, i.e. the restriction
	$\overline{j}_{|\pi (\mathfrak{C})}$ has a bijective differential at
	$\lambda $.
	
	On the other hand if $\pi (\mathfrak{C})$ and $\overline{j}^{-1}(z)$ are
	tangent at $\lambda $, then $d \overline{j}$ vanishes on
	$T_{\lambda}(\pi (\mathfrak{C}))$. This second case never occurs, since
	we have concluded above that
	$\overline{j}_{|\pi (\mathfrak{C}) \setminus S}: \pi (\mathfrak{C})
	\setminus S \rightarrow \mathbb C\setminus \{0,1728\}$ is a submersion.
	
	Therefore for any $z \neq 0,1728$, the fiber $\overline{j}^{-1}(z)$ and
	$\pi (\mathfrak{C})$ intersect transversely and at each intersection point
	the restriction $\overline{j}_{|\pi (\mathfrak{C}) \setminus S}$ has a
	bijective differential. Thus at each intersection
	$\overline{j}_{|\pi (\mathfrak{C}) \setminus S}$ defines a local biholomorphic
	bijection and there are $4$ such intersection points in $\mathcal{M}$ (or equivalently $24$ points in $\mathbb{C} \mathbb{P}^2$). So
	$\overline{j}_{|\pi (\mathfrak{C}) \setminus S}$ indeed defines a holomorphic
	covering space of degree $4$ over $\mathbb C\setminus \{0,1728\}$.
\end{proof}
	
From this, we can lift the result to $j$ itself, yielding the following
theorem.

%t5 #&#
\begin{theorem}
	%%LEAP%%%\label{thm5}
	\label{thm::fiber_bundle}
	Let $\mathfrak{S} = j^{-1}(0) \cup j^{-1}(1728) = \pi ^{-1}(S)$. Then the
	restriction (that we still denote $j$ for simplicity)
	$j: \mathfrak{C} \setminus \mathfrak{S} \rightarrow \mathbb C
	\setminus \{0,1728\}$ is a surjective submersion. More precisely, outside
	$\mathfrak{S}$, $j$ defines a fiber bundle over
	$\mathbb C\setminus \{0,1728\}$, the general fiber of which is isomorphic
	to the disjoint union $\coprod _{i=1}^{4} G/V_{4}$.
\end{theorem}
\begin{proof}
	First the surjectivity of $j_{|\mathfrak{C} \setminus \mathfrak{S}}$ needs
	to be proven. By Proposition~\ref{prop::covering},
	$\overline{j}_{|\pi (\mathfrak{C}) \setminus S}$ is surjective. Then
	$j = \overline{j} \circ \pi $ is a surjection from
	$\pi ^{-1}(\pi (\mathfrak{C}) \setminus S) = \mathfrak{C} \setminus
	\mathfrak{S}$ to $\mathbb C\setminus \{0,1728\}$.
	
	For $z \neq 0,1728$, we have
	$j^{-1}(z) \cap \mathfrak{C} = \pi ^{-1}(\overline{j}^{-1}(z) \cap
	\pi (\mathfrak{C}))$. Therefore by Proposition~\ref{prop::covering},
	$j^{-1}(z) \cap \mathfrak{C}$ has $4$ connected components, each being
	isomorphic to $G/V_{4}$ by Theorem~\ref{thm::orbits}, since the special
	orbit is in the fiber over $0$. Each of these components is the fiber of
	$\pi $ over an intersection point of $\overline{j}^{-1}(z)$ with
	$\pi (\mathfrak{C})$. Therefore the fibers of $j$ are all biholomorphic
	to $\coprod _{i=1}^{4} G/V_{4}$.
	
	Let $U$ be an open neighborhood of $z \neq 0,1728$, that is evenly covered
	through $\overline{j}$. Considering a smaller neighborhood if necessary,
	we can assume that $U$ is an open disk with center $z$. We are going to
	prove that
	$j_{|\mathfrak{C} \setminus \mathfrak{S}}^{-1}(U) \cong U \times
	\coprod _{i=1}^{4} G/V_{4}$.
	
	Let
	$\overline{j}_{|\pi (\mathfrak{C}) \setminus S}^{-1}(z) = \{p_{1},
	\ldots , p_{4}\}$ and let $V_{i}$ be the open neighborhood of
	$p_{i}$ in $\pi (\mathfrak{C}) \setminus S$ on which
	$\overline{j}_{|\pi (\mathfrak{C}) \setminus S}$ is biholomorphic onto
	$U$.
	
	We shall need the following notation. For $z'$ a point in $U$, we shall
	consider the segment $[z,z']$ parametrized by
	$\delta _{z'}: [0,1] \rightarrow U$, such that
	$\delta _{z'}(0) = z'$ and $\delta _{z'}(1) = z$.
	
	Let $q_{1} = p_{i} = [\lambda _{1}:\lambda _{2}:\lambda _{3}]$ and
	$q_{0} = [\mu _{1}:\mu _{2}:\mu _{3}]$ both in $V_{i}$, such that
	$j(q_{1}) = z$ and $j(q_{0}) = z'$. Let
	$\gamma _{z'}: [0,1] \rightarrow V_{i}$, the lift of $\delta _{z'}$, that
	is $\gamma _{z'} = \overline{j}_{|V_{i}}^{-1} \circ \delta _{z'}$. Consider
	the pairs of conics
	$\overline{q}_{1} = (I_{3},\text{diag}(\lambda _{1},\lambda _{2},
	\lambda _{3}))$ and
	$\overline{q}_{0} = (I_{3},\text{diag}(\mu _{1},\mu _{2},\mu _{3}))$. We
	have: $\pi (\overline{q}_{i}) = q_{i}$. Consider the path
	$\tilde{\gamma}_{z'}(t) = (I_{3},\text{diag}(\gamma _{z'}(t)))$. Therefore
	if $A \in \mathbb PGL_{3}(\mathbb C)$, let
	$\tilde{\gamma}_{z'}(t).A$ be
	$(A^{T}A,A^{T}\text{diag}(\gamma _{z'}(t))A)$.
	
	This allows to define a biholomorphic map:
	\begin{equation*}
		\begin{aligned}
			&j^{-1}(U) \rightarrow U \times j^{-1}(z),
			\\
			&(C,D) \mapsto (z',(A^{T} A , A^{T}\text{diag}(\lambda _{1},\lambda _{2},
			\lambda _{3})A) = (\overline{j}\circ \pi (C,D), \tilde{\gamma}_{z'}(1).A),
		\end{aligned}
	\end{equation*}
	where
	$(C,D) = (A^{T} A , A^{T} \text{diag}(\mu _{1},\mu _{2},\mu _{3}) A)$ and
	$\overline{j}([\mu _{1}: \mu _{2}: \mu _{3}]) = z'$.
	
	Therefore the function $j$ defines over $U$ a trivial fiber bundle, whose
	general fiber is $j^{-1}(z)$.
	
	Finally the inverse image $j^{-1}(U)$ is biholomorphic to
	$\bigcup _{i=1}^{4} (V_{i} \times G/V_{4}) \simeq U \times \coprod _{i=1}^{4}
	G/V_{4}$, which shows that $j$ indeed define a fiber bundle on
	$j: \mathfrak{C} \setminus \mathfrak{S} \rightarrow \mathbb C
	\setminus \{0,1728\}$, which generic fiber is
	$\coprod _{i=1}^{4} G/V_{4}$.
\end{proof}
	
Now we shall investigate what can be learned on the topology of
$\mathfrak{C}$ from this theorem. More precisely, we want to give some
information about the fundamental group of
$\mathfrak{C} \setminus \mathfrak{S}$.
	
\begin{theorem}
	\label{thm6}
	The fundamental group of $\mathfrak{C} \setminus \mathfrak{S}$ is the semidirect
	product $(V_{4} \times \mathbb Z_{3}) \rtimes _{\theta} K$, $K$ is a subgroup
	of $\mathbb{F}_{2}$, and
	$\theta : K \rightarrow \mathrm{Aut}(V_{4} \times \mathbb Z_{3})$ is a
	group homomorphism. Here $\mathbb{F}_{2}$ is the free group generated by
	$2$ elements.
\end{theorem}
\begin{proof}
	By Theorem~\ref{thm::fiber_bundle}, we will consider the following fiber
	bundle:
	\begin{equation*}
		F = \coprod _{i=1}^{4} G/V_{4} \hookrightarrow E = \mathfrak{C}
		\setminus \mathfrak{S} \overset{j}{\rightarrow} B = \mathbb C
		\setminus \{0,1728\} \rightarrow 0.
	\end{equation*}
	
	Let us consider $b_{0} \in B$ and $x_{0} \in j^{-1}(b_{0})$.
	
	First remind that by Proposition~\ref{prop::fund_group_orbits}, the fundamental
	group of $G = \mathbb PGL_{3}(\mathbb C)$ is
	$\mathbb Z_{3} = \mathbb Z/3\mathbb Z$, and the fundamental group of
	$G/V_{4}$ is $V_{4} \times \mathbb Z_{3}$
	
	Now consider the fiber $F = \coprod _{i=1}^{n} G/V_{4}$, with $n=4$. It
	has $n$ connected components, each being biholomorphic to $G/V_{4}$. Then
	for any base point $x_{0} \in F$, we have
	$\pi _{1}(F,x_{0}) = V_{4} \times \mathbb Z_{3}$.
	
	The base $B$ is homotopic to $S^{1} \vee S^{1}$. Therefore we have
	$\pi _{1}(B,x_{0}) = \mathbb Z\star \mathbb Z= \mathbb{F}_{2}$, the free
	group generated by two elements. Also, we have: $\pi _{2}(B) = 0$.
	
	Now we shall use again the long homotopy sequence, as in~\cite[th. 4.41]{hatcher-02}:
	\begin{equation*}
		\ldots \rightarrow \pi _{2}(B) \rightarrow \pi _{1}(F) \rightarrow
		\pi _{1}(E) \rightarrow \pi _{1}(B) \rightarrow \pi _{0}(F)
		\rightarrow \pi _{0}(E) \rightarrow 0,
	\end{equation*}
	since $B$ is connected. Note that by Lemma~\ref{lem::connectivity_codim2},
	$E$ is also connected, because it is the complement of
	$\mathfrak{S}$ in $\mathfrak{C}$ and the real codimension of
	$\mathfrak{S}$ in $\mathfrak{C}$ is $2$.
	
	Since $\pi _{2}(B) = 0$, we finally have the following exact sequence:
	\begin{equation*}
		1 \rightarrow V_{4} \times \mathbb Z_{3} \overset{i}{\rightarrow}
		\pi _{1}(E) \overset{\phi}{\rightarrow} \mathbb{F}_{2}
		\overset{\psi}{\rightarrow} \pi _{0}(F) \rightarrow 1
	\end{equation*}
	
	Let
	$K = \text{im}(\phi ) \cong \pi _{1}(E) / (V_{4} \times \mathbb Z_{3})$.
	Therefore we have:
	%
	%e10 #&#
	\begin{equation}
		\label{exact_seq_K}
		1 \rightarrow V_{4} \times \mathbb Z_{3} \overset{i}{\rightarrow}
		\pi _{1}(E) \overset{\phi}{\rightarrow} K \rightarrow 1
		%%LEAP%%%\label{eq10}
	\end{equation}
	
	Finally since $K$ is free, the sequence~\eqref{exact_seq_K} splits and
	this yields
	$\pi _{1} (\mathfrak{C} \setminus \mathfrak{S}) \cong (V_{4} \times
	\mathbb Z_{3}) \rtimes _{\theta} K$, with
	$\theta : K \rightarrow \mathrm{Aut}(V_{4} \times \mathbb Z_{3})$ being
	some group homomorphism.
\end{proof}
	
\section{Poncelet porism and J-invariant}
\label{sec::poncelet_principal_bundle}

For two non-degenerate conics that meet transversely, consider the Poncelet
correspondence:
\begin{equation*}
	\mathfrak{M}_{C,D} = \{(p,l) \in C \times D^{*} | p \in l\},
\end{equation*}
where $D^{*}$ denotes the dual conic of $D$. As proven in
\cite{griffiths_harris_78} and more explicitly in
\cite{gonzalez-16}, this set is actually an elliptic curve in
$\mathbb C\mathbb P^{1} \times \mathbb C\mathbb P^{1}$.

To prove this result, one typically proceeds as follows. First one shows
that $\mathfrak{M}_{C,D}$ is an algebraic curve of
$\mathbb C\mathbb P^{1} \times \mathbb C\mathbb P^{1}$ by using the parametric
representations of $C$ and $D^{*}$. The assumption that the conics meet
transversely ensures that this curve is non-singular. Finally considering
either projection
$\pi _{1}: C \times D^{*} \rightarrow \mathbb C\mathbb P^{1}, (p,l)
\mapsto p$ or
$\pi _{2}: C \times D^{*} \rightarrow \mathbb C\mathbb P^{1}, (p,l)
\mapsto l$, we have a holomorphic covering
$\mathcal{M}_{C,D} \rightarrow \mathbb C\mathbb P^{1}$ of degree $2$ with
four ramification points given by the intersection points of the two conics.
Then by Hurwitz formula, this yields that the genus of
$\mathcal{M}_{C,D}$ is 1. All together, the Poncelet correspondence
$\mathcal{M}_{C,D}$ is indeed an elliptic curve.

More concretely, one can observe that it is isomorphic as a Riemann surface
to
$E_{C,D} = \{(r,u) \in \mathbb C\mathbb P^{1} \times \mathbb C
\mathbb P^{1} | u_{0}^{2} r_{1}^{3} = u_{1}^{2} \det (r_{0}C+r_{1} D)
\}$.
	
It is therefore a genus one Riemann surface and it is isomorphic to a complex
torus. In particular, it is a complex Lie group. As a consequence there
exists a lattice $\Lambda $ such that
$\mathfrak{M}_{C,D} \simeq E_{C,D}$ is isomorphic to
$\mathbb C/\Lambda $. We shall denote this isomorphism by
$\nu _{C,D}: \mathfrak{M}_{C,D} \simeq \mathbb C/\Lambda $. Notice that
$\nu _{C,D}$ is not only a holomorphic function, it is also a group homomorphism,
provided the neutral element of $\mathfrak{M}_{C,D}$ is chosen to be
$\theta _{C,D} = \nu _{C,D}^{-1}([0])$, where $[0]$ is the class of zero
in $\mathbb C/\Lambda $, which we shall always assume in the sequel.
	
For a given value $z \neq 0,1728$ of the j-invariant, the set of isomorphic
elliptic curves, which satisfy the Cayley condition when regarded as points
in $\mathfrak{O}$, is given by
$j_{|\mathfrak{C} \setminus \mathfrak{S}}^{-1}(z)$, which has $4$ connected
components by Theorem~\ref{thm::fiber_bundle}. Consider the component
$\omega $ that contains $E_{C,D}$, which is nothing else than the orbit
of $(C,D)$ by the action of $\mathbb PGL_{3}(\mathbb C)$. Therefore
$\omega \simeq \mathbb PGL_{3}(\mathbb C)/ V_{4}$.

Let $\mathfrak{M}$ let be disjoint union
$\mathfrak{M} = \coprod _{(C,D) \in \omega} \mathfrak{M}_{C,D}$. Consider
the projection:
\begin{equation*}
	\tau : \mathfrak{M} \rightarrow \omega , (p,l) \mapsto (C,D),
\end{equation*}
if $(p,l) \in \mathfrak{M}_{C,D}$.
	
\begin{theorem}
	\label{thm7}
	The projection $\tau : \mathfrak{M} \rightarrow \omega $ defines a trivial
	principal fiber bundle, which general fiber is the complex torus
	$\mathbb C/\Lambda $.
\end{theorem}
\begin{proof}
	Define the function
	$\nu : \mathfrak{M} \rightarrow \omega \times (\mathbb C/\Lambda ), (p,l)
	\mapsto (\tau (p,l) , \nu _{\tau (p,l)}(p,l))$. This is a bijection which
	inverse is given by $\nu ^{-1}((C,D),[z]) = \nu _{C,D}^{-1}([z])$, where
	$[z]$ is the class of $z \in \mathbb C$ modulo $\Lambda $.
	
	This bijection endows $\mathfrak{M}$ with the structure of a complex manifold
	which makes it isomorphic to $\omega \times (\mathbb C/\Lambda )$. Since
	we have: $\tau = \text{pr}_{1} \circ \nu $ (where
	$\text{pr}_{1}: \omega \times (\mathbb C/\Lambda ) \rightarrow \omega $
	is the canonical projection) and $\nu _{C,D}$ is also a group homomorphism
	for all $(C,D) \in \omega $, this yields the result.
\end{proof}

\textit{Remark:} The triviality of this bundle is also reflected by the
existence of a global smooth section:
$s:\omega \rightarrow \mathfrak{M}, (C,D) \mapsto \theta _{C,D}$. This
section is holomorphic, since $s = \nu ^{-1} \circ i_{1}$, where
$i_{1}: \omega \rightarrow \omega \times (\mathbb C/\Lambda ), (C,D)
\mapsto ((C,D),[0])$ is the canonical injection.

\section{Conclusion}
\label{sec8}

In this paper, we have established that the regular Cayley set for triangles
is a smooth 9-dimensional complex manifold and identified its fiber bundle
structure over the punctured $j$-line. This investigation provides a rigorous
differential-geometric foundation for understanding the spaces of conic
pairs in Poncelet's porism.

Future work will be devoted to extending this differential-topological
approach to general values of $n$, seeking to characterize the global manifold
structure and topological invariants for $n$-gons beyond the case of triangles.
Additionally, the singular configurations that define the boundary of the
regular Cayley manifold will be investigated with the same differential-topological
flavor, with the aim of understanding how the manifold's topology degenerates
as it approaches these singular limits.

\medskip
{\bf Acknowledgment - } I would like to thank Prof. Misha Katz from Bar-Ilan University, for introducing me to the beautiful subject of Poncelet Porisms and for fruitful discussions. Also, I am thankful to Prof. Pierre Lochak from Institut Math\'ematique de Jussieu, Sorbonne Universit\'e for sharing with me his helpful commentaries.

%----------------------------------------------------------------------------------------
%	BIBLIOGRAPHY
%----------------------------------------------------------------------------------------
\bibliographystyle{numeric}

\end{document}